\documentclass[12pt]{article}

\usepackage{amsmath}
\usepackage{amsthm}
\usepackage{amssymb}
\usepackage[left=2.5cm,right=2.4cm,top=1.5cm,bottom=1.8cm]{geometry}
\usepackage{tikz}
\usetikzlibrary{positioning,shapes,calc,backgrounds,fit,matrix}
\usepackage{pgfplots}
\usetikzlibrary{plotmarks}
\usepackage{algpseudocode}
\usepackage{algorithm}
\usepackage{hyperref}

\numberwithin{equation}{section}
\numberwithin{table}{section}
\numberwithin{figure}{section}
\numberwithin{algorithm}{section}

\theoremstyle{definition}

\newtheorem{theorem}{Theorem}[section]
\newtheorem{lemma[theorem]}{Lemma}
\newtheorem{definition}[theorem]{Definition}
\newtheorem{remark}[theorem]{Remark}

\def\u{\mathbf{u}}
\def\A{\mathbf{A}}
\def\f{\mathbf{f}}
\def\c{\mathbf{c}}
\def\eps{\varepsilon}
\pgfsetlayers{background,main}

\pgfplotsset{%
  every axis/.append style={%
    y label style={at={(0.1,1.0)},anchor=south west,rotate=-90,color=black},
    yminorgrids
  }
}

\date{July 5, 2018}

\let\OLDthebibliography\thebibliography
\renewcommand\thebibliography[1]{
  \small
  \OLDthebibliography{#1}
  \setlength{\parskip}{0ex}
  \setlength{\itemsep}{0.5ex}
}

\begin{document}

\title{A hybrid Alternating Least Squares -- TT Cross algorithm for parametric PDEs\thanks{The first author acknowledges funding from the EPSRC fellowship EP/M019004/1.}}
\author{Sergey Dolgov and Robert Scheichl\footnote{University of Bath, Claverton Down, BA2 7AY, Bath, United Kingdom. {\tt s.dolgov@bath.ac.uk}, {\tt r.scheichl@bath.ac.uk}.}}

\maketitle

\begin{abstract}
We consider the approximate solution of parametric PDEs using the low-rank Tensor Train (TT) decomposition.
Such parametric PDEs arise for example in uncertainty quantification problems in engineering applications.
We propose an algorithm that is a hybrid of the alternating least squares and the TT cross methods.
It computes a TT approximation of the whole solution, which is beneficial when multiple quantities of interest are sought.
This might be needed, for example, for the computation of the probability density function (PDF) via the maximum entropy method [Kavehrad and Joseph, IEEE Trans. Comm., 1986].
The new algorithm exploits and preserves the block diagonal structure of the discretized operator in stochastic collocation schemes.
This disentangles computations of the spatial and parametric degrees of freedom in the TT representation.
In particular, it only requires solving independent PDEs at a few parameter values, thus allowing the use of existing high performance PDE solvers.
In our numerical experiments, we apply the new algorithm to the stochastic diffusion equation and compare it with preconditioned steepest descent in the TT format, as well as with (multilevel) quasi-Monte Carlo and dimension-adaptive sparse grids methods.
For sufficiently smooth random fields the new approach is orders of magnitude faster.

{\par\it Keywords:} stochastic PDEs, high-dimensional problems, tensor decompositions, low-rank decompositions, cross approximation.
\end{abstract}

\def\figscale{0.45}

\section{Introduction}
As a model, we consider the parameter-dependent diffusion problem
\begin{equation}
\begin{split}
 & -\nabla_x \cdot \Big( c(x,\mathbf{y}) \nabla_x u(x,\mathbf{y})\Big) = f(x,\mathbf{y}) \quad \mbox{with} \quad x \in \Omega \subset \mathbb{R}^{d_0},
\end{split}
\label{eq:pdeproblem}
\end{equation}
subject to appropriate boundary conditions,
where $d_0=2,3$ is the dimension of the physical space, and $\mathbf{y} = (y^1,\ldots,y^d)$ is a parameter, which can be high- or even infinite-dimensional.
A typical application is the modelling of stochastic subsurface flow, where the parametrization arises through a Karhunen-Lo\`eve expansion of the stochastic coefficient:
\begin{equation}
c(x,\mathbf{y}) = c_0(x) + \sum_{k=1}^{d} y^k \psi_k(x), \qquad y^k \sim \mathcal{U}(-\sqrt{3},\sqrt{3}).
\label{eq:kle_aff}
\end{equation}
In practice, of greater interest are the so-called log-normal and log-uniform cases,
\begin{equation}
c(x,\mathbf{y}) = \exp\left(\sum_{k=1}^{d} y^k \psi_k(x)\right),
\label{eq:kle}
\end{equation}
where $y^k \sim \mathcal{N}(0,1)$ or $y^k \sim \mathcal{U}(-\sqrt{3},\sqrt{3})$ are independent, normally or uniformly distributed random variables, respectively.
Typically, we are interested in evaluating some statistics of the solution $u$ itself, or of a function of $u$ -- a {\em quantity of interest (QoI)}.
To carry out the quadrature w.r.t.~the parameters requires solving problem \eqref{eq:pdeproblem} for a sufficiently large number of parameter values $\mathbf{y}$.
The parameter dimension $d$ can range from tens to thousands.
Classical approaches that rely on tensorizing one-dimensional quadrature rules are unthinkable in that case, since the number of quadrature points would grow exponentially in $d$.

The most common alternatives are (quasi) Monte Carlo \cite{nieder-qmc-1978,graham-QMC-2011} and quadrature methods based on sparse grid/polynomial approximation \cite{griebel-sparsegrids-2004}. Classical
Monte Carlo quadrature introduces randomly distributed points in a $d$-dimensional space, and converges independently of~$d$ (if the variance of the QoI is bounded independently of $d$).
However, the error decays no faster than $\mathcal{N}^{-1/2}$ for $\mathcal{N}$ points, which can be prohibitively slow.
Quasi-Monte Carlo (QMC) methods introduce a deterministic set of points, and can have higher convergence rates.
Standard QMC rules provide an error close to $\mathcal{O}(\mathcal{N}^{-1})$, and
this rate can be proved rigorously with a constant that is again independent of $d$ if the coefficient $c(x,\mathbf{y})$ satisfies certain conditions (both for the coefficients in \eqref{eq:kle_aff} and in \eqref{eq:kle}).
There exist higher order QMC rules \cite{Schwab-HOQMC-2014} which, under even stronger assumptions, can be shown to achieve faster convergence rates (but so far only for \eqref{eq:kle_aff}).
Another approach (with rigorous theory for \eqref{eq:kle_aff} and \eqref{eq:kle}) is the multilevel QMC method \cite{Schwab-MLQMC-2015,Scheichl-mlqmc-lognorm-2017}, which relies on a hierarchy of spatial discretizations: the (fast) solution on a coarse spatial mesh is used as a control variate for the solution on a finer mesh,
and hence fewer expensive fine-mesh solves are needed.

Quadrature methods based on sparse polynomial approximations for the solution of stochastic PDEs are usually constructed using stochastic collocation \cite{babuska-collocation-2007,nobile-sg-mc-2016} or stochastic Galerkin \cite{Bieri2009,matthies-galerkin-2005} techniques.
In particular, total-degree approximation builds on tensorized univariate polynomial interpolation with a maximum polynomial degree of $n$ in each coordinate direction.
In its basic form, it can only be used for moderate-dimensional problems (with $d$ up to $10$), since the total-degree polynomial space has a cardinality of $\mathcal{O}(\frac{(d+n)!}{d!n!})$, which is still very large for high $d$ and a typical size of $n$ of about $10$.
However, in the PDE application considered here, the different components of $\mathbf{y}$ typically have decreasing influence, since $\|\psi_k\|$ typically decreases algebraically or even exponentially with $k \to \infty$. In that case, adaptive sparse grids can be designed that require significantly smaller numbers of quadrature points \cite{Gerstner-adapt-SG-2003,spinterp}, and that can even be shown to converge dimension-independently when $c(x,\mathbf{y})$ is of the form \eqref{eq:kle_aff} in special cases \cite{schwab-estimates-2010}. No such proofs exist for the coefficients in \eqref{eq:kle}.

Tensor product decompositions form another class of adaptive techniques.
They trace back to separation of variables and the Fourier method for solving PDEs.
The solution can be sought in the form of products of univariate factors, e.g. $u(x,\mathbf{y}) \approx \sum_{\alpha=1}^{r} v_{\alpha}(x) w_{\alpha}(\mathbf{y})$,
where the number of terms $r$ is as small as possible.
Similarly, different components of $\mathbf{y}$ can be separated.
If $r$ is moderate, it becomes possible to use the tensor product quadrature and its superior accuracy,
since a high-dimensional integral breaks into one-dimensional integrals over each factor.
In practical computations, all functions are discretized and the separation of variables results in the so-called \emph{low-rank tensor formats}, generalizations of low-rank matrices.
For extensive reviews on the topic see \cite{kolda-review-2009,hackbusch-2012,bokh-surv-2015}.
Tensor product algorithms employ robust tools of linear algebra, such as the singular value decomposition and pivoted Gaussian elimination to deliver an optimal low-rank approximation for a given accuracy.
Several approaches to apply them to the stochastic PDE problem exist already \cite{KhSch-Galerkin-SPDE-2011,dklm-tt-pce-2015,Ballani-HTUQOut-2015,Schneider-HT-SFEM-2016,Ballani-HTUQRb-2016}.

Tensor product methods have the computational complexity
$\mathcal{O}(dnr^p)$, which scales linearly in the dimension $d$ and in the univariate grid size $n$, but polynomially in the \emph{ranks} $r$ that define the number of elements in tensor representation, and thereby the approximation accuracy.
In turn, the order $p$ and the hidden constant may be significantly different for different representations and algorithms.
In this paper, we use the Tensor Train (TT) decomposition \cite{osel-tt-2011}, which is both simple and robust, but generalisations, such as the hierarchical tensor (HT) format \cite{hk-ht-2009,gras-hsvd-2010}, exist and could also be employed.
The TT format has a storage complexity bounded by $dnr^2$, where $r$ is the largest TT rank.
However, the algorithmic complexity can feature $p>2$ and/or constants much greater than $1$. The overhead may be particularly high if the algorithm discards some of the underlying structure in the problem.
For example, \eqref{eq:pdeproblem} with finite element discretization in $x$ and collocation on a tensor product grid in $\mathbf{y}$ results in a large linear system with a block-diagonal structure (see Section \ref{sec:operator}).
Generic algorithms for the solution of tensorized linear systems destroy this structure in the course of the computation.
In that case, the computational complexity is $\mathcal{O}(d n^3r^6)$ for a direct solver or $\mathcal{O}(N_{i} dnr^4)$ for an iterative solver for the reduced linear systems (with hidden constants being in the order of 1), where $N_{i}$ is the number of iterations that can be in the order of thousands.
Moreover, if $N$ is the number of spatial degrees of freedom, the computation of the first TT factor (which carries the $x$ variable) requires the solution of a fairly dense $Nr \times Nr$ system, which does not lend itself easily to efficient preconditioning techniques.

We propose a new algorithm, where the TT factors corresponding to the $\mathbf{y}$ variables can be computed using a direct solver, but in a number of operations no bigger than $\mathcal{O}(dnr^4)$, with the hidden constant being again below 10.
The TT factor corresponding to the spatial variable $x$ can be computed using only $\mathcal{O}(r)$ independent deterministic PDE solves and $\mathcal{O}(Nr^3)$ additional floating point operations.
Extensive numerical experiments have shown that the TT rank $r$ grows only logarithmically with the required accuracy $\eps$ (in the approximation of the QoI) and at most linearly in $d$.
The actual value of $r$ depends on the smoothness of the input random field $c(x,\mathbf{y})$ -- our experiments include fairly rough models, similar to lognormal fields with Mat\'ern covariance with smoothness parameter $\nu=1$. Even using efficient preconditioning techniques, the $\mathcal{O}(|\log \eps|)$ PDE solves will typically dominate the computational cost, especially in three spatial dimensions. Thus, the overall cost estimate is much lower than for QMC methods, which require $\mathcal{N}=\mathcal{O}(\eps^{-1})$ PDE solves to achieve an accuracy of $\eps$.

Our overall TT scheme is carried out as follows.
Firstly,  using the TT-Cross algorithm~\cite{ot-ttcross-2010}, we approximate the whole coefficient $c(x,\mathbf{y})$ at all grid points in $x$ and $\mathbf{y}$ in the TT representation (see Section \ref{sec:cross}).
The TT-Cross algorithm requires to evaluate $c(x,\mathbf{y})$ only at $\mathcal{O}(dnr^2)$ grid points in $x$ and $\mathbf{y}$,
where $r$ grows at most linearly in $d$ and logarithmically in~$\eps$.
As already mentioned, so far there is no rigorous proof for this statement but substantial numerical evidence.
Then, using our proposed hybrid method based on the TT-Cross and an Alternating Least Squares (ALS) linear solver \cite{holtz-ALS-DMRG-2012}, we approximate the whole PDE solution $u(x,\mathbf{y})$ at all the grid points in the TT representation (see Section \ref{sec:als}), providing a so-called \emph{surrogate model} or \emph{response surface}. Since a pointwise evaluation of a TT decomposition is inexpensive, the accuracy of a TT approximation can be certified cheaply using, for example, a Monte Carlo estimate. We will provide such evidence with our numerical experiments.
Finally, multiple statistics of multiple QoIs can be evaluated at negligible cost by direct tensor product integration based on the TT approximation of $u(x,\mathbf{y})$.

In \cite{Ballani-HTUQOut-2015}, the authors use a cross algorithm to directly build a surrogate model in HT format (in the high-dimensional parameter $\mathbf{y}$) for a functional of the PDE solution (e.g. $\int_{\Omega} u(x,\mathbf{y})dx$), avoiding the approximation of the entire solution $u$ at all the spatial grid points. When only one moment of one QoI is needed this is likely going to lead to a lower rank in the TT approximation and to a lower computational cost than our proposed TT scheme. However, often several statistics of several QoIs are required. One such example is the Maximum Entropy method \cite{Chernov-maxentr-2016} for approximating the entire probability density function (PDF) of a QoI~$Q$, which is useful for computing probabilities of events, including rare ones, or quantiles. The quality of the approximation, especially of the tails of the PDF, depends strongly on the number of moments.
It may thus be necessary to compute many moments, see Section \ref{sec:num_PDF}.
In that case, the TT rank of the simultaneous approximation of all those moments may be of the same order as (or even higher than) the rank required to approximate the whole solution. Moreover, once the whole solution is available in a structured representation, further QoIs or higher moments are available immediately at negligible extra cost.

The paper is organised as follows.
Section \ref{sec:tt} contains the definition and some basic properties of TT decomposition, as well as the TT-Cross and the Alternating Least Squares algorithm for the approximation and for the solution of linear systems in TT format, respectively.
In Section \ref{sec:alscross-main} we present a combination of the two algorithms and tailor it to provide an efficient algorithm for discretized parameter-dependent PDEs.
Finally, Section \ref{sec:numer_main} is devoted to numerical experiments, where we compare the newly proposed TT algorithm with single and multi-level QMC and Sparse Grid techniques.
We demonstrate that, for sufficiently fast decaying $\|\psi_k\|$ in the expansions \eqref{eq:kle_aff} and \eqref{eq:kle}, the new algorithm is consistently up to two orders of magnitude faster.

\section{Background on TT algorithms}
\label{sec:tt}

\subsection{Tensor Train decomposition and notation}
In this section, we give a very brief introduction to the tensor train methodology in an abstract notation.
To link with the original PDE problem in Section \ref{sec:alscross-main},
we deviate slightly from the common notation in tensor literature.
We consider a \emph{tensor} $\mathbf{v}(j_0,\ldots,j_d)$, where the total number of variables (the \emph{dimension}) is $d+1$,
and each $j_k$ for $k=0,\ldots,d$ enumerates degrees of freedom (DOFs) corresponding to the $k$-th variable.
On the other hand, $\mathbf{v}$ must also be represented as a vector to apply standard linear algebra tools.
An equivalence between tensors and vectors is established via multi-indices.
\begin{definition}
Given indices $j_0,\ldots,j_d$ with ranges $n_0,\ldots,n_d$,
 a \emph{multi-index} is introduced as their straightforward lexicographic grouping,
$$
\overline{j_0\ldots j_d} = (j_0-1) n_1 \cdots n_d + \cdots + (j_{d-1}-1)n_d + j_d.
$$
\end{definition}
Now we can represent the same data by either a tensor $\mathbf{v}(j_0,\ldots,j_d)$ or a vector $v(\overline{j_0\ldots j_d})$.

Throughout the paper, we use the Tensor Train (TT) decomposition \cite{osel-tt-2011} to represent (or approximate) tensors:
\begin{equation}
 \mathbf{v}(j_0,\ldots,j_d) = \sum_{\alpha_0,\ldots,\alpha_{d-1}=1}^{r_0,\ldots,r_{d-1}}\mathbf{v}^{(0)}_{\alpha_0}(j_0)\mathbf{v}^{(1)}_{\alpha_0,\alpha_1}(j_1) \cdots \mathbf{v}^{(d)}_{\alpha_{d-1}}(j_d).
\label{eq:tt}
\end{equation}
Each element of $\mathbf{v}$ is represented (or approximated) by a sum of products of elements of smaller tensors $\mathbf{v}^{(k)} \in \mathbb{R}^{r_{k-1} \times n_k \times r_k}$, called \emph{TT blocks}.
The auxiliary summation indices $\alpha_0,\ldots,\alpha_{d-1}$ are called rank indices, and their ranges $r_0,\ldots,r_{d-1}$ are called \emph{TT ranks}.
For simplicity, we define $r_{-1} = r_{d} = 1$.
In the physics literature, the TT decomposition is known as the Matrix Product States \cite{schollwock-2005} representation, since in \eqref{eq:tt},
each element of $\mathbf{v}$ is given by a product of $d+1$ matrices, with the $k$-th matrix depending on the ``state'' of the index $j_k$.
The TT ranks depend on the approximability of the particular tensor and on the required accuracy, if \eqref{eq:tt} does not hold exactly.
Defining upper bounds $r_k \le r$, $j_k \le n$, $k=0,\ldots,d$, we can see that the
required storage for the TT blocks is $\le dnr^2 + nr$, which can be much less than the full amount $n^{d+1}$.
We refer to $r$ as the maximal TT rank.

If a tensor is given fully or as another TT representation, a TT approximation with quasi-optimal TT ranks for a given accuracy threshold $\eps$ can be computed via the Singular Value Decomposition \cite{osel-tt-2011}.
This is crucial to avoid growth of the TT ranks when performing
algebraic operations, such as additions or multiplications, with TT representations.

A matrix $A \in \mathbb{R}^{n^{d+1} \times n^{d+1}}$ can be seen as a $(2d+2)$-dimensional tensor and represented in a slightly different TT format,
\begin{equation}
A(\overline{j_0\ldots j_d},\overline{j_0' \ldots j_d'}) = \sum\limits_{\gamma_0,\ldots,\gamma_{d-1}=1}^{R_0,\ldots,R_{d-1}} \A^{(0)}_{\gamma_0}(j_0,j_0') \A^{(1)}_{\gamma_0,\gamma_1}(j_1,j_1') \cdots  \A^{(d)}_{\gamma_{d-1}}(j_d,j_d'),
\label{eq:ttm}
\end{equation}
where $\A^{(k)} \in \mathbb{R}^{R_{k-1} \times n_k \times n_k \times R_k}$ are the \emph{matrix} TT blocks.
This is consistent with the Kronecker product ($\otimes$) in two dimensions.
The matrix-vector product $Av$ can then be implemented in the TT format block by block \cite{osel-tt-2011}.

In what follows, to shorten calculations it will be convenient to consider a TT block $\mathbf{v}^{(k)}$ (a 3-dimensional tensor) as a matrix or a vector.
\begin{definition}
 Given a $3$-dimensional tensor (TT block) $\mathbf{v}^{(k)} \in \mathbb{R}^{r_{k-1} \times n_k \times r_k}$,
 we define the following repartitions of its elements:
 \begin{enumerate}
  \item[(a)] \emph{Vector} folding: $\mathrm{v}^{(k)}(\overline{\alpha_{k-1} j_k \alpha_k}) = \mathbf{v}^{(k)}_{\alpha_{k-1},\alpha_k}(j_k)$, \quad $\mathrm{v}^{(k)} \in \mathbb{R}^{r_{k-1}n_kr_k \times 1}$.
  \item[(b)] \emph{Left} matrix folding:  $V^{|k\rangle}(\overline{\alpha_{k-1} j_k}, \alpha_k) = \mathbf{v}^{(k)}_{\alpha_{k-1},\alpha_k}(j_k)$, \quad $V^{|k\rangle} \in \mathbb{R}^{r_{k-1}n_k \times r_k}$.
  \item[(c)] \emph{Right} matrix folding: $V^{\langle k|}(\alpha_{k-1}, \overline{j_k \alpha_k}) = \mathbf{v}^{(k)}_{\alpha_{k-1},\alpha_k}(j_k)$, \quad $V^{\langle k |} \in \mathbb{R}^{r_{k-1}\times n_k r_k}$.
 \end{enumerate}
Equivalently, once any of the foldings is defined, we assume that the original tensor, as well as all others foldings are also defined and carry the same elements.
\label{def:folds}
\end{definition}

The full range (a vector of all values) of an index $i$ is denoted by $\{i\}$ or ``$:$'', for example, $A(:,j)=A(\{i\},j)$ is the $j$-th column of a matrix $A$ with elements $A(i,j)$.

\subsection{TT-Cross approximation}
\label{sec:cross}

Even if a low-rank TT approximation exists, it is not always obvious how to obtain it efficiently.
In this section, we recall a general approach for building TT
approximations that requires only a small number of evaluations of tensor elements.
We assume that we do not have access to the whole tensor, but only a
procedure that can return elements of the tensor at any arbitrary index.

The TT-Cross algorithm \cite{ot-ttcross-2010} is an extension of Gaussian elimination to TT tensors.
An $n \times m$ matrix $V$ of rank $r$ can be recovered by computing a
so-called \emph{cross}
\begin{equation}
V = V(:,\mathcal{J}) V(\mathcal{I},\mathcal{J})^{-1} V(\mathcal{I},:)
\label{eq:matrixcross}
\end{equation}
of $r$ columns and rows, where $\mathcal{I} \subset \{1,\ldots,n\}$
and $\mathcal{J} \subset \{1,\ldots,m\}$ are two index sets of
cardinality $r$ such that $V(\mathcal{I},\mathcal{J})$ (the intersection matrix) is nonsingular.
If $r\ll n$ and $m$, this decomposition requires computing only $(n+m-r)r \ll nm$ elements of the original matrix.
If the matrix is not exactly of low rank, we can still try to approximate it with the right hand side of \eqref{eq:matrixcross}.
However, we need a stable algorithm for computing such $\mathcal{I}$ and $\mathcal{J}$ that not only provide a nonsingular intersection matrix, but minimize also the approximation error.

A theoretically optimal, but NP-hard choice of
$\mathcal{I}$ and $\mathcal{J}$  is
such that $V(\mathcal{I},\mathcal{J})$ has
the maximum \emph{volume} (= modulus of the determinant) among all $r \times r$ submatrices.
Different sub-optimal strategies have been suggested in practice, such as ACA \cite{bebe-2000},
a search restricted to a random submatrix \cite{Ballani-HTUQOut-2015} and alternating maximum volume iteration \cite{gostz-maxvol-2010}.
In this paper we employ the latter approach, which is particularly convenient for combining with other tensor product methods (see the next section).

The idea is to replace the (hard) problem of finding a $r \times r$
optimal submatrix in a $n \times m$ matrix by a feasible problem of
finding a $r \times r$ submatrix in a $n \times r$ matrix, when $r \ll m$. Let one
index set be given; without loss of generality, let it be $\mathcal{J}$.
Then it is feasible to compute the $nr$ elements of the corresponding $r$ columns $V(:,\mathcal{J})$.
Now a $r \times r$ submatrix that has close to maximum volume can be found in
$V(:,\mathcal{J})$ instead of in the full matrix $V$ by the so-called
\emph{maxvol} algorithm \cite{gostz-maxvol-2010} in $\mathcal{O}(nr^2)$ operations.
We collect the corresponding row positions of this submatrix into an index set
$$
\mathcal{I} = \mathrm{maxvol}(V(:,\mathcal{J})) = \left\{i_{\alpha}\right\}_{\alpha=1}^{r}, \quad \mbox{s.t.} \quad |\det V(\mathcal{I},\mathcal{J})| \ge c \max\limits_{\mathcal{\hat I} \subset \{i\},\#\mathcal{\hat I}=r} |\det V(\mathcal{\hat I},\mathcal{J})|, \quad c>0.
$$
The index set $\mathcal{I}$ can in turn be used to compute $r$ rows $V(\mathcal{I},:)$.
The maxvol algorithm, applied to $V(\mathcal{I},:)^\top$, gives a \emph{new} index set $\mathcal{J}$,
and the process repeats.
Since the volume of $V(\mathcal{I},\mathcal{J})$ is non-decreasing and bounded from above, it converges.
A few of such \emph{alternating} iterations are usually sufficient to obtain index sets that deliver good approximations.

However, especially in the first iterations, the factor $V(:,\mathcal{J})$ can be ill-conditioned or even rank-deficient.
This makes it difficult to determine the maximum volume index set, not to mention that the inversion in \eqref{eq:matrixcross} becomes impossible.
Fortunately, it is sufficient to orthogonalize $V(:,\mathcal{J})$ to stabilize the calculations.
Indeed, in exact arithmetic,
$$
V(:,\mathcal{J}) V(\mathcal{I},\mathcal{J})^{-1} = Q R (Q(\mathcal{I},:) R)^{-1} = QRR^{-1} Q(\mathcal{I},:)^{-1} = QQ(\mathcal{I},:)^{-1},
$$
where $QR = V(:,\mathcal{J})$ is the QR decomposition, and hence $\mathrm{maxvol}(Q)$ and $\mathrm{maxvol}(V(:,\mathcal{J}))$ produce exactly the same $\mathcal{I}$.
But in the presence of round-off errors, using the $Q$-factor gives a much more robust algorithm in comparison to \eqref{eq:matrixcross}.

A full TT decomposition can now be computed by applying this alternating maxvol algorithm recurrently.
To start the procedure, suppose we are given an initial \emph{right} index set $\mathcal{J}_{>0} \in \mathbb{N}^{r_0 \times d}$, containing rows of indices $(j_1,\ldots,j_d)_{\alpha_0}$, for $\alpha_0=1,\ldots,r_0$.
Then it is feasible to compute $n_0 r_0$ elements $\mathbf{v}(j_0, (j_1)_{\alpha_0}, \ldots, (j_d)_{\alpha_0})$, to collect them into a matrix $V^{|0\rangle}$, to orthogonalize it and to apply the maxvol procedure to deduce the \emph{left} index set $\mathcal{I}_{<1}$.
Now, suppose in the $k$-th step we are given index sets $\mathcal{I}_{<k} \in \mathbb{N}^{r_{k-1} \times k}$ and $\mathcal{J}_{>k} \in \mathbb{N}^{r_k \times d-k}$, containing some instances of $(j_0,\ldots,j_{k-1})_{\alpha_{k-1}}$ and $(j_{k+1},\ldots,j_d)_{\alpha_k}$, respectively.
We can evaluate $r_{k-1} n_k r_k$ elements
\begin{equation}
V^{| k\rangle }(\overline{\alpha_{k-1} j_k}, \alpha_k) = \mathbf{v}((j_0)_{\alpha_{k-1}}, \ldots, (j_{k-1})_{\alpha_{k-1}}, j_k, (j_{k+1})_{\alpha_k} \ldots, (j_d)_{\alpha_k}),
\label{eq:newttblock}
\end{equation}
orthogonalize this matrix and compute the set of \emph{local} maxvol indices
$$
\mathcal{L}_k = \mathrm{maxvol}(V^{| k\rangle }) = \left\{(\overline{\alpha_{k-1} j_k})_{\alpha_k}\right\}_{\alpha_k=1}^{r_k}.
$$
Since $V^{| k\rangle }$ consists of the elements of $\mathbf{v}$,
$\alpha_{k-1}$ is associated with the tuples $(j_0,\ldots,j_{k-1})_{\alpha_{k-1}}$,
and we can map $(\alpha_{k-1})_{\alpha_k}$ from $\mathcal{L}_k$ to the set of \emph{global} indices
\begin{equation}
\mathcal{I}_{<k+1} = \left\{(j_0)_{(\alpha_{k-1})_{\alpha_k}}, \ldots, (j_{k-1})_{(\alpha_{k-1})_{\alpha_k}}, (j_k)_{\alpha_k} \right\}_{\alpha_k=1}^{r_k}, \quad \left(j_0,\ldots,j_{k-1}\right)_{\alpha_{k-1}} \in \mathcal{I}_{<k}.
\label{eq:indexmerge-l}
\end{equation}
This procedure continues until $k=d$, and then recurs in a similar fashion \emph{backward} from $k=d$ to $k=0$, updating the right indices in the process:
\begin{equation}
\mathcal{J}_{>k-1} = \left\{ (j_k)_{\alpha_{k-1}}, (j_{k+1})_{(\alpha_k)_{\alpha_{k-1}}} \ldots, (j_d)_{(\alpha_k)_{\alpha_{k-1}}} \right\}_{\alpha_{k-1}=1}^{r_{k-1}}, \quad \left(j_{k+1},\ldots,j_d\right)_{\alpha_k} \in \mathcal{J}_{>k},
\label{eq:indexmerge-r}
\end{equation}
where $(j_k)_{\alpha_{k-1}}$ and $(\alpha_k)_{\alpha_{k-1}}$ belong to
a similar set of local maxvol indices
$$
\mathrm{maxvol}((V^{\langle k| })^\top) = \left\{(\overline{j_k \alpha_{k}})_{\alpha_{k-1}}\right\}_{\alpha_{k-1}=1}^{r_{k-1}}.
$$

Instead of specifying initial indices $\mathcal{J}_{>k}$, we can start with any TT decomposition as an initial guess.
In this case, we compute $\mathcal{J}_{>k}$ from the maxvol procedure applied to the initial TT blocks in the first iteration.
The whole procedure is outlined in Algorithm \ref{alg:ttcross}.

\begin{algorithm}[t]
\caption{TT-Cross algorithm \cite{ot-ttcross-2010}}
\label{alg:ttcross}
\begin{algorithmic}[1]
\Require Initial TT blocks $\mathbf{v}^{(k)}$, $k=0,\ldots,d;$ function to evaluate $\mathbf{v}(j_0,\ldots,j_d)$; tolerance $\eps$.
\Ensure Updated TT approximation blocks $\mathbf{v}^{(k)}$.

\State Initialize $\mathcal{I}_{<0}=\mathcal{J}_{>d}=\emptyset$ (see also Remark \ref{rem:rand-cross}) and $\mbox{iter}=0$.
\While{$\mbox{iter}<I_{\max}$ or $\|\mathbf{v} - \mathbf{v}_{prev}\| > \eps \|\mathbf{v}\|$}
  \State Set $\mbox{iter} \leftarrow \mbox{iter}+1$ and $\mathbf{v}_{prev} \leftarrow  \mathbf{v}$
  \For{$k=0,1,\ldots,d$}
   \If{$\mathcal{J}_{>k}$ exist, e.g. $\mbox{iter}>1$}
    \State\label{al:cross:evl} Evaluate $\mathbf{v}^{(k)} = \mathbf{v}\left(\mathcal{I}_{<k} \cup \{j_k\} \cup \mathcal{J}_{>k}\right) \in \mathbb{R}^{r_{k-1} \times n_k \times r_k}$, as given in \eqref{eq:newttblock}.
   \EndIf
    \If{$k<d$}
      \State\label{al:cross:qr} Compute QR decomposition $V^{| k \rangle} = QR$, \quad $Q \in \mathbb{R}^{r_{k-1}n_k \times r_k}$.
      \State\label{al:cross:mvl} Determine a new vector of local pivots $\mathcal{L}_{k} = \mathrm{maxvol}(Q)$.
      \State\label{al:cross:Vl} Update TT blocks: $V^{| k \rangle} \leftarrow QQ(\mathcal{L}_k,:)^{-1}$ and $V^{\langle k+1|} \leftarrow Q(\mathcal{L}_k,:) R V^{\langle k+1|}$
      \State\label{al:cross:Il} Map $\mathcal{L}_k$ to global index set $\mathcal{I}_{<k+1} = \left[\mathcal{I}_{<k} \cup \{j_k\}\right]_{\mathcal{L}_k}$, as given in \eqref{eq:indexmerge-l}.
    \EndIf
  \EndFor
  \For{$k=d,d-1,\ldots,0$}
      \State\label{al:cross:evr} Evaluate $\mathbf{v}^{(k)} = \mathbf{v}\left(\mathcal{I}_{<k} \cup \{j_k\} \cup \mathcal{J}_{>k}\right) \in \mathbb{R}^{r_{k-1} \times n_k \times r_k}$, as given in \eqref{eq:newttblock}.
    \If{$k>0$}
      \State\label{al:cross:lq} Compute QR decomposition $(V^{\langle k|})^\top = QR$, \quad $Q \in \mathbb{R}^{n_k r_k \times r_{k-1}}$.
      \State\label{al:cross:mvr} Determine a new vector of local pivots $\mathcal{L}_{k} = \mathrm{maxvol}(Q)$
      \State\label{al:cross:Vr} Update TT blocks: $V^{\langle k |} \leftarrow Q(\mathcal{L}_k,:)^{-\top} Q^\top$ and $V^{| k-1 \rangle} \leftarrow V^{| k-1 \rangle} R^\top Q(\mathcal{L}_k,:)^{\top}$.
      \State\label{al:cross:Jr} Map $\mathcal{L}_k$ to the global index set $\mathcal{J}_{>k-1} = \left[\{j_k\} \cup \mathcal{J}_{>k}\right]_{\mathcal{L}_{k}}$, as given in \eqref{eq:indexmerge-r}.
    \EndIf
  \EndFor
\EndWhile
\end{algorithmic}
\end{algorithm}

\begin{remark}
\label{rem:rand-cross}
In general, the selection of the initial guess, or, equivalently, of the initial indices, might be a crucial, but difficult task.
For example, consider a tensor of all zeros except for one single nonzero entry.
The probability of finding its location (an initial index) at random is $1/n^d$, i.e. extremely small.
However, this is to some extent an academic example.
The tensors in our application are smooth.
In this case, if no specific initial guess is available,
a reasonable choice is to initialize $\mathcal{J}_{>k}$ with $\mathcal{N}$ random indices,
where $\mathcal{N}$ is slightly larger than the expected ranks.
We can therefore evaluate Line~\ref{al:cross:evl} already in the first iteration.
To compress the solution to optimal ranks,
we replace the (exact) QR decomposition in Line~\ref{al:cross:qr} by an approximate truncated singular value or cross decomposition.
The random index initialization was also employed in the cross approximation in the HT format \cite{Ballani-HTUQOut-2015,Ballani-HTUQRb-2016}.
However, the HT cross is implemented differently: both $\mathcal{I}_{<k}$ and $\mathcal{J}_{>k}$ are sought among a priori chosen random subsets of the corresponding spaces.
In Alg.~\ref{alg:ttcross}, one of the sets ($\mathcal{I}_{<k}$) is optimized by the maxvol algorithm (in particular, $\mathcal{I}_{<1}$ is chosen from the \emph{entire} range $\{i_0\}$).
Moreover, if the initial representation is insufficient, better index sets are found by conducting several alternating iterations.
\end{remark}

\subsection{Alternating Least Squares linear solver}
\label{sec:als}

Cross algorithms are suitable when an explicit function for evaluating an arbitrary element of a tensor is available.
However, this is not the case if the tensor is given as the unknown solution of a linear system, such as $A\mathrm{u}=\mathrm{f}$,
where $\mathrm{u}$ is a $n^{d+1} \times 1$ vectorization of the sought tensor  $\mathbf{u}(j_0,\ldots,j_d)$.
A central element of the TT-Cross algorithm is the alternating iteration over the TT blocks.
A similar idea applies also to linear systems, leading to the so-called Alternating Least Squares (ALS) algorithm \cite{holtz-ALS-DMRG-2012}.

We assume that $A$ is symmetric positive definite and that it can be represented in TT format as in \eqref{eq:ttm}.
Let $\mathrm{f}$ also be represented in TT-format with ranks $\rho_0,\ldots,\rho_d$, i.e.,
\begin{equation*}
\begin{split}
 \mathrm{f}(\overline{j_0\ldots j_d}) & = \sum_{\delta_{0},\ldots,\delta_{d-1}=1}^{\rho_0,\ldots,\rho_{d-1}} \f^{(0)}_{\delta_0}(j_0) \f^{(1)}_{\delta_0,\delta_1}(j_1) \cdots \f^{(d)}_{\delta_{d-1}}(j_d).
\end{split}
\end{equation*}

We reformulate the solution of the linear system $A\mathrm{u}=\mathrm{f}$ as a minimization problem for the functional $J(\mathrm{u}) = \mathrm{u}^\top A\mathrm{u} -2\mathrm{u}^\top \mathrm{f}$ and plug the TT decomposition of $\mathrm{u}$
into the functional~$J$.
We solve this minimisation problem by iterating over $k=0,\ldots,d$ and restricting the optimization in each step to a single TT block:
\begin{equation}
\u^{(k)} = \arg\min_{\mathbf{v}^{(k)} \in \mathbb{R}^{r_{k-1} \times n_k \times r_k}} J(\mathrm{u}), \quad \mbox{where}
\label{eq:localopt}
\end{equation}
\begin{equation}
\mathrm{u}(\overline{j_0\ldots j_d}) = \u^{(0)}(j_0) \cdots \u^{(k-1)}(j_{k-1}) \mathbf{v}^{(k)}(j_k)  \u^{(k+1)}(j_{k+1}) \cdots \u^{(d)}(j_d).
\label{eq:tt_u}
\end{equation}
In order to solve the univariate minimisation problem \eqref{eq:localopt} for each $k$, we rewrite \eqref{eq:tt_u} as a single matrix product.
\begin{definition}
\label{def:iface}
By \emph{interface matrices}, or simply \emph{interfaces}, we denote partial TT decompositions, with elements\vspace{-1ex}
\begin{equation}
\begin{split}
 U^{(<k)}_{\alpha_{k-1}}(\overline{j_0\ldots j_{k-1}}) & = \sum_{\alpha_0,\ldots,\alpha_{k-2}=1}^{r_0,\ldots,r_{k-2}} \u^{(0)}_{\alpha_0}(j_0) \cdots \u^{(k-1)}_{\alpha_{k-2},\alpha_{k-1}}(j_{k-1}), \quad U^{(<k)} \in \mathbb{R}^{n_0\cdots n_{k-1} \times r_{k-1}}, \\
 U^{(>k)}_{\alpha_{k}}(\overline{j_{k+1}\ldots j_d}) & = \sum_{\alpha_{k+1},\ldots,\alpha_{d-1}=1}^{r_{k+1},\ldots,r_{d-1}} \u^{(k+1)}_{\alpha_k,\alpha_{k+1}}(j_{k+1}) \cdots \u^{(d)}_{\alpha_{d-1}}(j_d), \quad U^{(>k)} \in \mathbb{R}^{r_{k} \times n_{k+1}\cdots n_d}.
\end{split}
\label{eq:iface}
\end{equation}
To simplify the presentation, we also define $U^{(<0)}=1$ and $U^{(>d)}=1$.
\end{definition}

Using these interfaces, \eqref{eq:tt_u} can be represented as
$$
\mathrm{u} = U_{\neq k} \mathrm{v}^{(k)}, \quad \text{where} \ \ U_{\neq k} = U^{(<k)} \otimes I_{n_k} \otimes (U^{(>k)})^\top,
$$
$I_{n_k}$ is the identity matrix of size $n_k$, and $\mathrm{v}^{(k)}$ is a $r_{k-1}n_kr_k \times 1$ vectorization of the TT block $\mathbf{v}^{(k)}$ (see Def. \ref{def:folds}).
The so-called \emph{frame} matrix $U_{\neq k} \in \mathbb{R}^{n^{d+1} \times r_{k-1}n_k r_k}$ defines thus a linear map from the $k$th TT block to the whole tensor.
Finally, we can rewrite \eqref{eq:localopt} as a reduced symmetric positive definite $r_{k-1}n_kr_k \times r_{k-1}n_kr_k$ linear system, i.e.,
\begin{equation}
\mathrm{u}^{(k)} = \arg\min_{\mathrm{v}^{(k)}\in\mathbb{R}^{r_{k-1}n_kr_k}} (U_{\neq k} \mathrm{v}^{(k)})^\top AU_{\neq k} \mathrm{v}^{(k)} - 2 (U_{\neq k} \mathrm{v}^{(k)})^\top \mathrm{f} \ \ \Leftrightarrow \ \ \left(U_{\neq k}^\top AU_{\neq k}\right) \mathrm{u}^{(k)} = U_{\neq k}^\top \mathrm{f}.
\label{eq:localsys}
\end{equation}

The representations of $A$, $\mathrm{f}$ and $U_{\neq k}$ in TT format can be exploited for an efficient assembly of the reduced system \eqref{eq:localsys}.
In particular, $\mathrm{f}_k = U_{\neq k}^\top \mathrm{f}$ reduces to
\begin{equation}
 \mathrm{f}_k = \sum_{\delta_{k-1},\delta_k=1}^{\rho_{k-1},\rho_k}\hat F_{<k}(:,\delta_{k-1}) \otimes \mathbf{f}^{(k)}_{\delta_{k-1},\delta_k}(:) \otimes \hat F_{>k}(\delta_k, :) \in \mathbb{R}^{r_{k-1}n_k r_k}, \quad \mbox{where} \ \ \hat F_{<k} = (U^{(<k)})^\top F^{(<k)}
 \label{eq:flocal}
\end{equation}
and  $\hat F_{>k} = F^{(>k)} (U^{(>k)})^\top$
can be computed by exploiting the TT formats of $\mathrm{u}$ and~$\mathrm{f}$ again.
These projections are the counterparts of the index sets $\mathcal{I}_{<k}$ and $\mathcal{J}_{>k}$ in the TT-Cross Algorithm (Alg.~\ref{alg:ttcross}).
The reduced matrix can be constructed in a similar way:
\begin{equation}
 A_k = U_{\neq k}^\top AU_{\neq k} = \sum_{\gamma_{k-1},\gamma_k=1}^{R_{k-1},R_k}\hat A_{<k}(:,:,\gamma_{k-1}) \otimes \mathbf{A}^{(k)}_{\gamma_{k-1},\gamma_k} \otimes \hat A_{>k}(\gamma_k,:,:) \in \mathbb{R}^{r_{k-1}n_k r_k \times r_{k-1}n_k r_k},
 \label{eq:Alocal}
\end{equation}
where $\mathbf{A}^{(k)}_{\gamma_{k-1},\gamma_k} = \mathbf{A}^{(k)}_{\gamma_{k-1},\gamma_k}(:,:)$ is a matrix slice of the TT block $\mathbf{A}^{(k)}$.
To define the projections $\hat A_{<k}$, $\hat A_{>k}$ in \eqref{eq:Alocal},
we introduce the interfaces of $A$, which are defined (similarly to those of vectors in Def. \ref{def:iface}) as follows:
\begin{equation}
A^{(<k)}_{\gamma_{k-1}} = \sum_{\gamma_0,\ldots,\gamma_{k-2}} \A^{(0)}_{\gamma_0} \otimes  \cdots \otimes \A^{(k-1)}_{\gamma_{k-2},\gamma_{k-1}} \ \ \text{and} \ \
A^{(>k)}_{\gamma_k} = \sum_{\gamma_{k+1},\ldots,\gamma_{d-1}} \A^{(k+1)}_{\gamma_k,\gamma_{k+1}} \otimes  \cdots \otimes \A^{(d)}_{\gamma_{d-1}}.
\label{eq:iface-mat}
\end{equation}
The projections $\hat A_{<k}$, $\hat A_{>k}$ are then simply Galerkin projections of
these interfaces onto the solution interfaces, i.e.,
\begin{equation}
\hat A_{<k}(:,:,\gamma_{k-1}) = (U^{(<k)})^\top A^{(<k)}_{\gamma_{k-1}} U^{(<k)} \ \ \text{and} \ \  \hat A_{>k}(\gamma_{k},:,:) = U^{(>k)} A^{(>k)}_{\gamma_{k}} (U^{(>k)})^\top.
\label{eq:Aleft}
\end{equation}
In the course of the alternating iteration, these projections can be computed recursively, using only the $k$-th TT blocks of $A$, $\mathrm{f}$ and $\mathrm{u}$ and the projections from the previous step \cite{holtz-ALS-DMRG-2012}.

However, the solution interfaces $U^{(<k)}$, $U^{(>k)}$ are dense.
Hence, the projections in \eqref{eq:Aleft} and the reduced matrix in \eqref{eq:Alocal}
are also dense, even if the TT blocks $\A^{(k)}$ are sparse, e.g. even if $\A^{(k)}(j_k,j_k')=0$ for all $j_k \neq j_k'$.
The next section is devoted to the main contribution of this paper,
on how to avoid this unnecessary ``fill-in'' and the associated computational cost.

\section{ALS-Cross algorithm for parameter-dependent PDEs}
\label{sec:alscross-main}

In this section we introduce our new hybrid algorithm, where the iteration for $k=0,1,\ldots,d$ is carried out in the ALS fashion, and the backward iteration for $k=d,\ldots,0$ resembles the backward iteration of the TT-Cross algorithm.
However, first we need to discretize the parametric PDE \eqref{eq:pdeproblem} both in space and in the parameters.

\subsection{Discretization of the elliptic parametric PDE}
For the spatial discretization we use continuous piecewise linear (or bilinear) finite elements.
We denote by $\mathcal{T}_h$ a shape-regular family of simplicial triangulations or rectangular tessellations of the Lipschitz polygonal/polyhedral
domain $\Omega$, parametrized by the mesh size $h = \max_{\tau \in \mathcal{T}_h} \mathrm{diam}(\tau)>0$.
The partitioning $\mathcal{T}_h$ induces a finite element (FE) space
$
\bar V_h = \text{span}\left\{\phi_s(x): s=1,\ldots,\bar N \right\}
$
of continuous
functions on $\Omega$ that are linear (resp. bilinear) on each element
$\tau \in \mathcal{T}_h$.
For every fixed $\mathbf{y} \in \mathbb{R}^d$,
the FE approximation to the solution of \eqref{eq:pdeproblem} is the unique function $u_h \in V_h \subset \bar V_h$, satisfying the weak formulation
$$
\int_{\Omega} c(x,\mathbf{y}) \nabla\phi_i(x)\cdot \nabla u_h(x,\mathbf{y}) dx  = \int_{\Omega} f(x,\mathbf{y})\phi_i(x) dx, \quad \forall i=1,\ldots,N,
$$
where $V_h = \text{span}\left\{\phi_i(x): i=1,\ldots,N \right\}$ is the subspace of $\bar V_h$ containing only those FE functions that also satisfy the essential boundary conditions on $\partial\Omega$ (for details see Sec.~\ref{sec:numer_main}).
Existence and uniqueness of $u_h$, as well as the bound $|u_h(\cdot,\mathbf{y})|_{H^1(\Omega)} \le \|f(\cdot,\mathbf{y})\|_{H^{-1}(\Omega)}/\min_{x}c(x,\mathbf{y})$ follow by
the Lax-Milgram Lemma.
Standard convergence results are of the form
\[
|u - u_h|_{H^r(\Omega)} \le C h^{s-r}, \quad \text{where} \ \ r \in [0,1], \ \ s \in (0,2] \ \text{and} \ \ C>0
\]
depend on the smoothness of the coefficient, the boundary and the right hand side \cite{Hackbusch-elPDE-1992}. For $c\in W^{1,\infty}(\Omega)$ and $f \in L^2(\Omega)$ on a convex domain $\Omega$, we have $\|u-u_h\|_{L^2(\Omega)} \le C h^2$.

Since a closed form of the coefficient $c(x,\mathbf{y})$ is rarely available, we need to discretize it as well. For simplicity, in this paper, we interpolate $c(x,\mathbf{y})$ in the same FE space $\bar V_h$, such that $c_h(x,\mathbf{y}) = \sum_{s=1}^{\bar N} c(x_s,\mathbf{y}) \phi_s(x)$.
Similarly we discretize the right hand side.
Now the finite element solution is defined by $A(\mathbf{y}) \mathrm{u}(\mathbf{y}) = \mathrm{f}(\mathbf{y})$, where
\begin{equation}
A_{i,i'}(\mathbf{y}) = \int_{\Omega} \sum_{s=1}^{\bar N} c(x_s,\mathbf{y}) \phi_s(x) \nabla\phi_i(x)\cdot \nabla\phi_{i'}(x) dx, \quad \mathrm{f}_{i}(\mathbf{y}) = \int_{\Omega} \sum_{s=1}^{\bar N} f(x_s,\mathbf{y}) \phi_s(x) \phi_i(x) dx,
\label{eq:A_f_fix}
\end{equation}
for $i,i'=1,\ldots,N$, and $\mathrm{u}$ is a coefficient vector, representing $u_h \in V_h$.

The dependency on the parameters $\mathbf{y}$ is discretized via the stochastic collocation scheme on a tensor product grid.
If $\mathbf{y}$ is normally (resp.~uniformly) distributed, let $Y_k = \{y^k_{j_k}\}_{j_k=1}^{n_k} \subset (-a,a)$ denote the roots of the Hermite (resp.~Legendre) polynomials for each of the components $y^k$ of $\mathbf{y}$, where $a = \infty$ (resp.~$a=\sqrt{3}$). Using the sets $Y_k$, we construct the following set of $d$-dimensional grid points,
$$
\mathcal{Y}_n = \left\{\mathbf{y}_{\mathbf{j}}=(y^1_{j_1},\ldots,y^d_{j_d}) : y^k_{j_k} \in Y_k \quad \forall j_k=1,\ldots,n_k, \quad 1 \le n_k \le n, \quad k=1,\ldots,d \right\},
$$
which can be used to define the space
\begin{equation}
\mathcal{V}_n = \left\{L_{\mathbf{j}}(\mathbf{y}) = \prod_{k=1}^{d} L_{j_k}(y^k) : L_{j_k}(y^k) \mbox{ is a Lagrange polynomial centered at } y^k_{j_k} \in Y_k  \right\}
\label{eq:Vn}
\end{equation}
of globally continuous Lagrange interpolation polynomials on $(-a,a)^d$.
The approximate solution is then sought in the form $u_{h,n}(x,\mathbf{y}) = \sum_{\mathbf{y}_\mathbf{j} \in \mathcal{Y}_n} u_{\mathbf{j}}(x) L_{\mathbf{j}}(\mathbf{y})$ with $(u_{\mathbf{j}})_{\mathbf{y}_\mathbf{j} \in \mathcal{Y}_n} \subset V_h$.

Having discretized both in space and in the parameters, the discrete approximations of the coefficient, the right-hand side and the solution can be represented as the following \emph{tensors}:
$$
\mathbf{c}(s,j_1,\ldots,j_d) = c(x_s, \mathbf{y}_{\mathbf{j}}), \ \ \mathbf{f}(i, j_1,\ldots,j_d) = \mathrm{f}_i(\mathbf{y}_{\mathbf{j}}) \ \ \text{and} \ \  \u(i, j_1,\ldots,j_d) = u_{\mathbf{j}}(x_i) = \mathrm{u}_i(\mathbf{y}_{\mathbf{j}}),
$$
for $i=1,\ldots,N,~s=1,\ldots,\bar N,~j_k=1,\ldots,n_k,~k=1,\ldots,d$.

The $n_1n_2\ldots n_d$ unknown FE functions $(u_{\mathbf{j}})_{\mathbf{y}_\mathbf{j} \in \mathcal{Y}_n} \subset V_h$ are then uniquely defined by assuming that the semidiscrete FE system $A(\mathbf{y}) \mathrm{u}(\mathbf{y}) = \mathrm{f}(\mathbf{y})$ is satisfied at all the grid points $\mathbf{y}_\mathbf{j} \in \mathcal{Y}_n$ (collocation). This finally reduces the continuous problem \eqref{eq:pdeproblem} to a linear system with a block-diagonal matrix,
\begin{equation}
\underbrace{\begin{bmatrix}
 A(\mathbf{y}_{1,\ldots,1}) \\
 & \ddots \\
 & & A(\mathbf{y}_{n_1,\ldots,n_d})
\end{bmatrix}}_{A}
\begin{bmatrix}
 \u(:,1,\ldots,1) \\
 \vdots \\
 \u(:,n_1,\ldots,n_d)
\end{bmatrix}=
\begin{bmatrix}
 \mathbf{f}(:,1,\ldots,1) \\
 \vdots \\
 \mathbf{f}(:,n_1,\ldots,n_d)
\end{bmatrix},
\label{eq:block-diag}
\end{equation}
where $\mathbf{y}_{j_1,\ldots,j_d} \in \mathcal{Y}_n$, ``:'' denotes the full range of $i$ from $1$ to $N$, and $A$ and $\mathbf{f}$ are assembled according to \eqref{eq:A_f_fix}.
The block-diagonal form arises from the property of the Lagrange polynomials that $L_{\mathbf{j}}(\mathbf{y}_{\mathbf{k}}) = 0$, for all $\mathbf{j} \not=\mathbf{k}$.

However, even to store the tensors $\mathbf{c},\mathbf{u}$ or $\mathbf{f}$ directly would be prohibitive, since the amount of values $N \prod_k n_k \le Nn^d$ exceeds the memory capacity already for moderate $N,n$ and $d$.
Instead, we approximate  $\mathbf{c},\mathbf{u}$ and $\mathbf{f}$ using TT decompositions.

\subsection{Construction of the block-diagonal system in TT format}
\label{sec:operator}

First, the discrete coefficient $\c$ needs to be provided in TT format
\begin{equation}
\c(s,\mathbf{j}) = \sum_{\gamma_0,\ldots,\gamma_{d-1}=1}^{R_0,\ldots,R_{d-1}}\c^{(0)}_{\gamma_0}(s)\c^{(1)}_{\gamma_0,\gamma_1}(j_1) \cdots \c^{(d)}_{\gamma_{d-1}}(j_d), \quad s=1,\ldots,\bar N,~j_k=1,\ldots,n_k,~k=1,\ldots,d.
\label{eq:tt-kappa}
\end{equation}
For the affine form \eqref{eq:kle_aff}, the following natural rank-$d$ TT decomposition can be used:
\begin{equation}
\mathbf{c}(s,\mathbf{j}) = \sum_{k=1}^{d} \psi_k(x_s) \cdot e_{j_1} \cdots e_{j_{k-1}} \cdot y^k_{j_k} \cdot e_{j_{k+1}} \cdots e_{j_d},
\label{eq:c_aff}
\end{equation}
where $e$ is a vector of all ones.
We construct this decomposition directly.

For the non-affine coefficient \eqref{eq:kle} no such explicit analytical TT decomposition exists. Instead, we use the TT-Cross algorithm (Alg.~\ref{alg:ttcross}) to compute an approximate TT decomposition, passing the function $\mathbf{v}(s,j_1,\ldots,j_d) = \c(s,\mathbf{j}) = \exp\left(\sum_{k=1}^{d} \psi_k(x_s) y^k_{j_k}\right)$ for the evaluation of the elements of the discrete tensor $\c$ to the TT-Cross algorithm.
The initial indices for TT-Cross can be chosen randomly, see Remark \ref{rem:rand-cross}.

Having decomposed the coefficient, we still need to decompose the matrix and the right-hand side to solve the system \eqref{eq:block-diag}.
A crucial ingredient is the fact that the bilinear form in \eqref{eq:A_f_fix} depends linearly on the coefficient.
Plugging the TT format of the coefficient \eqref{eq:tt-kappa} into the bilinear form \eqref{eq:A_f_fix}, we obtain
\begin{equation}
A_{i,i'}(\mathbf{j}) = \left[\int_{\Omega} \sum_{s=1}^{\bar N} \c^{(0)}(s) \phi_s(x) \nabla\phi_i(x)\cdot \nabla\phi_{i'}(x) dx \right] \cdot  \c^{(1)}(j_1) \cdots \c^{(d)}(j_d).
\label{eq:tt-vm}
\end{equation}
These matrices $A(\mathbf{j})$ are the blocks on the diagonal of \eqref{eq:block-diag}.
Now, notice that for any matrix $P$ and vector $\mathrm{q} = (q_1,\ldots,q_n)$ we have
\begin{equation}
 \mathrm{diag}(\mathrm{q}) \otimes  P = \begin{bmatrix}P q_1 \\ & \ddots \\ & &  P q_n\end{bmatrix}, \quad \mbox{where} \quad \mathrm{diag}(\mathrm{q})_{j,j'} = q_j\delta_{j,j'},
\label{eq:kron-diag}
\end{equation}
and $\delta_{j,j'}$ is the Kronecker symbol, equal to $1$ when $j=j'$ and $0$ otherwise.
Due to distributivity, \eqref{eq:kron-diag} can be extended to the TT format \eqref{eq:tt-vm}.
Hence, the whole matrix $A$ in \eqref{eq:block-diag} admits a TT decomposition \eqref{eq:ttm} with the ranks $R_0,\ldots,R_{d-1}$ as in \eqref{eq:tt-kappa} and the following TT blocks:
\begin{equation}
\begin{split}
\A^{(0)}_{\gamma_0}(i,i') & \ = \ \int_{\Omega} \sum_{s=1}^{\bar N}\c^{(0)}_{\gamma_0}(s)\phi_s(x) \nabla\phi_i(x)\cdot \nabla\phi_{i'}(x)dx, \\
\A^{(k)}_{\gamma_{k-1},\gamma_k}(j_k,j_k')  & \ = \ \c^{(k)}_{\gamma_{k-1},\gamma_k}(j_k) \delta_{j_k,j_k'}, \qquad\qquad k>0.
\end{split}
\label{eq:A_pde}
\end{equation}

The right-hand side is assembled or approximated in the TT format as well, i.e.
$$
\mathbf{f}(i,\mathbf{j}) = \sum_{\delta_0,\ldots,\delta_{d-1}=1}^{\rho_0,\ldots,\rho_{d-1}}\mathbf{f}^{(0)}_{\delta_0}(i) \mathbf{f}^{(1)}_{\delta_0,\delta_1}(j_1) \cdots \mathbf{f}^{(d)}_{\delta_{d-1}}(j_d).
$$
The simplest case is a deterministic right-hand side, $f=f(x)$, which yields a rank-1 TT decomposition
$$
\mathbf{f}(i,\mathbf{j}) = \left[\int_{\Omega} \sum_{s=1}^{\bar N} f(x_s) \phi_s(x) \phi_i(x) dx \right] e_{j_1} \cdots e_{j_d}.
$$
Another example would be the case of inhomogeneous deterministic Dirichlet boundary conditions.
In that case, the right hand side depends linearly on the elements of the matrix $A$ and thus again linearly on the coefficient $c$,
leading to a rank-$(R_0,\ldots,R_{d-1})$ TT format of $\f$.

Before we proceed to the solution algorithm, it is instructive to see how an interpolation $u_{h,n} \in V_h \otimes \mathcal{V}_n$ can be evaluated efficiently for arbitrary points $x \in \Omega$ and $\mathbf{y} \in (-a,a)^d$, provided a TT decomposition of the tensor $\mathbf{u}$ containing the nodal values of $u_h$ is available.
Indeed, the indices $j_k$ in \eqref{eq:Vn} are independent, and the summation over $\mathbf{j}$ can be distributed with the multiplications in the TT format,
\begin{equation}
u_{h,n}(x,\mathbf{y}) = \left[\sum_{i=1}^{N} \u^{(0)}(i) \phi_i(x)\right] \left[\sum_{j_1=1}^{n_1} \u^{(1)}(j_1) L_{j_1}(y^1)\right] \cdots \left[\sum_{j_d=1}^{n_d} \u^{(d)}(j_d) L_{j_d}(y^d)\right].
\label{eq:lagrange}
\end{equation}
Since the finite element basis functions $\phi_i(x)$ are local, and $L_{j_k}(y^k)$ are global (univariate) functions, evaluation of $u_{h,n}(x,\mathbf{y})$ at fixed $x,\mathbf{y}$ requires $\mathcal{O}(r+dnr^2)$ operations.
Similarly, the \emph{quadrature weights} for integrating the interpolant $u_{h,n}$ over $x,\mathbf{y}$ can be obtained by explicit integration of $\phi_i$ and $L_{j_k}$, and the overall quadrature complexity is $\mathcal{O}(Nr + dnr^2)$.

\begin{remark}
The presented approach is not limited to the discrete description.
We can consider continuous variables, such that TT blocks become matrix-valued functions \cite{osel-constr-2013}, e.g.
$$
c(x,\mathbf{y}) = \c^{(0)}(x) \c^{(1)}(y^1) \cdots \c^{(d)}(y^d),
$$
or tensors with one of the dimensions being infinite, e.g. $\c^{(k)}_{\alpha_{k-1},\alpha_k}(y^k)$.
Matrix multiplications involving this dimension (e.g. ALS projections \eqref{eq:Aleft}) should be replaced by integrals.
Products over rank indices (in e.g. QR decompositions) imply linear combinations of columns of the corresponding matrix foldings,
which are to be replaced by linear combinations of functions.
Computation of cross indices can be also performed at the continuous level, as determining one index means finding the maximal value of a function in modulus.
This is done in the Empirical Interpolation method (EIM), for example \cite{Maday-eim-2004}.
From the practical perspective, it may pave the way for an adaptive discretization in space and in the parameters, interleaved with the ALS/cross iterations over the TT structure \cite{bachmayr-sparse-or-lr-2017}.
\end{remark}

\subsection{A hybrid ALS-Cross algorithm for solving the block-diagonal collocation system}
Our goal is to replace the dense, reduced matrix \eqref{eq:Alocal} in the ALS method by a sparser matrix which
allows an efficient solution.
Recall from \eqref{eq:A_pde} that, for $k=1,\ldots,d$, the TT factors $\A^{(k)}(j_k,j_k')=0$ whenever $j_k \neq j_k'$.
This structure is generally destroyed by the projection onto the dense interface $U^{(>k)}$ of the solution in \eqref{eq:Alocal}.

In order to sparsify the factor $\hat A_{>k}$ in \eqref{eq:Alocal} we modify the ALS algorithm in two aspects.
Firstly, instead of using the dense interface $U^{(>k)}$ in the $k$-th step directly,
we use a similar trick to that used in Alg.~\ref{alg:ttcross} at Lines \ref{al:cross:lq}--\ref{al:cross:Jr}, i.e., we recursively apply the maxvol algorithm to find an index set $\mathcal{J}_{>k}$ and update the $k$-th TT block.
This is equivalent to replacing $U^{(>k)}$ with $\big[U^{(>k)}(:,\mathcal{J}_{>k})\big]^{-1} U^{(>k)}$.
This transformation ensures that $U^{(>k)}(:,\mathcal{J}_{>k}) = I_{r_k}$.
Secondly, we replace the symmetric Galerkin projection $U_{\neq k}^\top A U_{\neq k}$ in \eqref{eq:localsys} by a Petrov-Galerkin projection $W_{\neq k}^\top A U_{\neq k}$ with a specially chosen matrix $W_{\neq k} \in \mathbb{R}^{N n_1\cdots n_d \times r_{k-1} n_k r_k}$.
Since we are interested in sparsifying the right projection $\hat A_{>k}$, we simply set the columns of $U^{(>k)}$ corresponding to indices $\mathbf{j} \not\in \mathcal{J}_{>k}$ to zero, i.e.
$$
W_{\neq k} = U^{(<k)} \otimes I_{n_k} \otimes \left(E^{(>k)}\right)^\top,
$$
where $E^{(>k)}  \in \mathbb{R}^{r_k \times n_{k+1}\cdots n_{d}}$ is a submatrix of identity at the index set $\mathcal{J}_{>k}$, that is
$$
E^{(>k)}(:,\mathcal{J}_{>k}) = I_{r_{k}},
$$
while all other elements of $E^{(>k)}$ are zero.
For the projected matrix we now obtain:
\begin{equation}
W_{\neq k}^\top A U_{\neq k} =
\sum_{\gamma_{k-1},\gamma_k=1}^{R_{k-1},R_k}\hat A_{<k}(:,:,\gamma_{k-1}) \otimes \mathbf{A}^{(k)}_{\gamma_{k-1},\gamma_k} \otimes \left(E^{(>k)} A^{(>k)}_{\gamma_k} (U^{(>k)})^\top\right) \in \mathbb{R}^{r_{k-1}n_k r_k \times r_{k-1}n_k r_k},
\label{eq:PGmat}
\end{equation}
where $A^{(>k)}$ is as given in \eqref{eq:iface-mat}. Finally, we solve the Petrov-Galerkin linear system
\begin{equation}
(W_{\neq k}^\top A U_{\neq k})\mathrm{u}^{(k)} = (W_{\neq k}^\top \mathrm{f})
\label{eq:PG}
\end{equation}
to compute the next TT block $\u^{(k)}$.

The block-diagonality of the \emph{hybrid} projection \eqref{eq:PGmat} can be shown as follows. Without loss of generality, let us assume that
$\mathcal{J}_{>k} = \{1,\ldots,r_k\}$, and note that due to \eqref{eq:A_pde}
$A^{(>k)}_{\gamma_k}$ is diagonal. Then, we obtain
$$
E^{(>k)} A^{(>k)}_{\gamma_k} (U^{(>k)})^\top =
\begin{bmatrix}
 A^{(>k)}_{\gamma_k}(1,1)     & \cdots & 0 & 0  \\
 \vdots & \ddots & \vdots & 0 \\
 0 & \cdots & A^{(>k)}_{\gamma_k}(r_k,r_k) & 0
\end{bmatrix}
\underbrace{\begin{bmatrix}
 1 & \cdots & 0 \\
 \vdots & \ddots & \vdots \\
 0 & \cdots & 1 \\
 * & * & *
\end{bmatrix}}_{(U^{(>k)})^\top} =
A^{(>k)}_{\gamma_k}(\mathcal{J}_{>k}, \mathcal{J}_{>k}).
$$
Hence, both the middle and the right factor in \eqref{eq:PGmat}, for each $\gamma_{k-1},\gamma_k$, are diagonal, and so the matrix $W_{\neq k}^\top A U_{\neq k}$ is block-diagonal.
Let us consider the cases $k=0$ and $k>0$ separately.

\begin{description}
 \item[$\mathbf{k=0:}$] For uniformity of notation, let $\hat A_{<k}=1$.
Assume again that $\mathcal{J}_{>0}=\{1,\ldots,r_0\}$ and recall the definition of the zeroth TT block $\A^{(k)}_{\gamma_0}=\A^{(0)}_{\gamma_0}$ in \eqref{eq:A_pde}.
It follows from \eqref{eq:PGmat}, since $W_{\neq 0} = I_N \otimes (E^{(>0)})^\top = \begin{bmatrix}I_{Nr_0} \!& \!0\end{bmatrix}$, that
 $$
 W_{\neq k}^\top A U_{\neq k} = \sum_{\gamma_0} \A^{(0)}_{\gamma_0} \otimes A^{(>0)}_{\gamma_0}(\mathcal{J}_{>0},\mathcal{J}_{>0}) =
 \begin{bmatrix}A\left(\mathbf{y}_{1}\right) \\ & \ddots \\ & & A\left(\mathbf{y}_{r_0}\right)\end{bmatrix},
 $$
 where $A(\mathbf{y}_{\mathbf{j}})$ is the FE stiffness matrix at parameter value $\mathbf{y}_{\mathbf{j}}$, for $\mathbf{j} \in \mathcal{J}_{>0}$, according to \eqref{eq:A_f_fix}.
 Similarly, the right hand side $W_{\neq k}^\top \mathrm{f}$ is simply the concatenation of the vectors $\mathrm{f}(\mathbf{y}_1),\ldots,\mathrm{f}(\mathbf{y}_{r_0})$, as defined in \eqref{eq:A_f_fix}.

Therefore, the step $k=0$ is similar to the offline stage in \emph{Reduced Basis methods} \cite{Ballani-HTUQRb-2016}:
we solve independent deterministic PDEs at a finite set of $r_0$ parameter values to produce the snapshots $\mathrm{u}(\mathbf{y}_{\mathbf{j}})$,
$\mathbf{j} \in \mathcal{J}_{>0}$, which then form the columns of the TT block $\u^{(0)}$.
The difference is that the choice of the parameter values is not random, but optimised by the maxvol algorithm. As the algorithm converges, the indices in $\mathcal{J}_{>0}$ will converge to their maximum volume positions. For smooth functions, they should deliver better approximation than random (Monte Carlo) samples.

In order to proceed to the next step ($k=1$), we need to construct the partially reduced factor $\hat A_{<1}$ in \eqref{eq:Aleft}.
This involves orthogonalisation of $\u^{(0)}$ and projection of $\A^{(0)}$ and $\f^{(0)}$ onto this orthogonal basis, as in Reduced Basis Methods.
In this step it is possible to reduce the rank $r_0$, using truncated singular value decomposition, similarly to Proper Orthogonal Decomposition. For more details on adapting the ranks see Section \ref{sec:amen}.

 \item[$\mathbf{k>0:}$] Due to the block-diagonal form of the matrix \eqref{eq:PGmat},
 for the other TT blocks, corresponding to the discretization in $\mathbf{y}$, the solution of \eqref{eq:PG} decouples into $n_kr_k$ independent dense systems of size $r_{k-1} \times r_{k-1}$.
 These blocks are orthogonal projections of the FE stiffness matrices at some $\mathbf{y}$ points, and hence are symmetric positive definite.
 This is similar to the \emph{online stage} of Reduced Basis Methods.
 The difference is that this online computation is split into $d$ steps over individual TT blocks.
\end{description}

The pseudocode of the ALS-Cross method is shown in Alg.~\ref{alg:als-cross}.
Similarly to the TT-Cross and ALS methods, this algorithm is iterated forward, from $k=0$ to $k=d$, and back, from $k=d$ to $k=0$.
In the forward iteration where $k$ increases, we compute the orthogonal projection $\hat A_{<k}$ as in \eqref{eq:Aleft} for the ALS method.
In the backward iteration where $k$ decreases, we use the recursive maxvol algorithm and update the TT blocks as described in this section.

\begin{algorithm}[t]
\caption{ALS-Cross algorithm}
\label{alg:als-cross}
\begin{algorithmic}[1]
\Require TT blocks of the coefficient $\c^{(k)}$ and right-hand side $\mathbf{f}^{(k)}$; initial index sets $\mathcal{J}_{>k}$, $k=0,\ldots,d$; stopping tolerance $\eps$.
\Ensure TT blocks of the solution $\u^{(k)}$.

\State Initialize $\hat A_{<0} = \hat F_{<0}=1$, copy $\u^{(k)} \leftarrow \c^{(k)}$, $k=0,\ldots,d$, set $\mbox{iter}=0$.
\While{$\mbox{iter}<I_{\max}$ or $\|\u - \u_{prev}\|>\eps \|\u\|$}
  \State Set $\mbox{iter} \leftarrow \mbox{iter}+1$ and $\u_{prev} \leftarrow \u$.
  \State\label{al:ac:det} For all $\alpha_0=1,\ldots,r_0,$ solve $A(\mathbf{y}_{\alpha_0}) \u^{(0)}_{\alpha_0} = \mathrm{f}(\mathbf{y}_{\alpha_0})$, with $A(\mathbf{y}_{\alpha_0}),\mathrm{f}(\mathbf{y}_{\alpha_0})$ as given in \eqref{eq:A_f_fix}.
  \State\label{al:ac:qr0} Orthogonalize $\u^{(0)} = QR$ and update the TT blocks $\u^{(0)} \leftarrow  Q$, $U^{\langle 1|} \leftarrow  RU^{\langle 1|}$.
  \State\label{al:ac:V0} Compute projections $\hat A_{<1}$ and $\hat F_{<1}$ as shown in \eqref{eq:Aleft} and \eqref{eq:flocal}.
  \For{$k=1,2,\ldots,d$}
    \State\label{al:ac:evl} Solve the reduced problem \eqref{eq:PG} with the matrix from \eqref{eq:PGmat}.
    \If{$k<d$}
      \State\label{al:ac:qr} Orthogonalize $U^{| k \rangle} = QR$ and update $U^{| k \rangle} \leftarrow Q$, $U^{\langle k+1|} \leftarrow RU^{\langle k+1|}$.
      \State\label{al:ac:Al} Compute projections $\hat A_{<k+1}$ and $\hat F_{<k+1}$ as shown in \eqref{eq:Aleft} and \eqref{eq:flocal}.
    \EndIf
  \EndFor
  \For{$k=d,d-1,\ldots,1$}
    \State\label{al:ac:evr} Solve the reduced problem $(W_{\neq k}^\top A U_{\neq k}) \mathrm{u}^{(k)} = W_{\neq k}^\top \mathrm{f}$, as above.
    \State\label{al:ac:lq} Orthogonalize $(U^{\langle k |})^\top = QR$.
    \State\label{al:ac:mvr} Determine a new vector of local pivots $\mathcal{L}_{k} = \mathrm{maxvol}(Q)$.
    \State\label{al:ac:Vr} Update the TT blocks $U^{\langle k |} \leftarrow Q(\mathcal{L}_k,:)^{-\top} Q^\top$ and $U^{| k-1 \rangle} \leftarrow U^{| k-1 \rangle} R^\top Q(\mathcal{L}_k,:)^{\top}$
    \State\label{al:ac:Jr} Map $\mathcal{L}_k$ to the global index set $\mathcal{J}_{>k-1} = \left[\{j_k\} \cup \mathcal{J}_{>k}\right]_{\mathcal{L}_k}$, as given in \eqref{eq:indexmerge-r}.
  \EndFor
\EndWhile
\end{algorithmic}
\end{algorithm}

\subsection{Adaptation of TT ranks}
\label{sec:amen}
The algorithms listed so far work with TT approximations of fixed rank.
However, the final ranks of the solution are difficult to guess a priori.
A practical algorithm should be able to decrease or increase the ranks depending on a desired accuracy.
To decrease the ranks is easy: we just need to apply any rank-revealing technique (singular value decomposition or matrix cross) \cite{holtz-ALS-DMRG-2012} to an appropriate folding of the TT block, $U^{|k\rangle}$ or $U^{\langle k |}$, after it is computed from the reduced system (Lines \ref{al:ac:evl} and \ref{al:ac:evr} in Alg.~\ref{alg:als-cross}).

More interesting is the rank-increasing step.
There are two main approaches.
The Density Matrix Renormalization Group (DMRG) method suggested in the quantum physics community \cite{white-dmrg-1993} merges two neighbouring TT blocks into one larger block, $\u^{(k,k+1)}$, then the corresponding larger system is solved, and the variables are separated back.
In the latter step, the rank $r_k$ can be increased up to a threshold.
An immediate drawback is the need to solve a larger ``two-dimensional'' system in each step.
A less obvious caveat is the possible stagnation of convergence far
from the optimal solution, which can occur especially for high
dimensional problems.
A remedy for this problem is to pad the TT blocks of the solution with random entries \cite{DoOs-dmrg-solve-2011}.

It was observed later that using this expansion the ranks can actually be increased explicitly without merging TT blocks.
A more reliable and theoretically justified method, where the TT blocks of the solution are augmented by the TT blocks of the residual, is the Alternating Minimal Energy (AMEn) method proposed in \cite{ds-amen-2014}.
It converges faster than DMRG and its complexity remains of the same order as that of the plain vanilla ALS algorithm.

The AMEn algorithm features a \emph{semi-reduced} system $A_{\ge k} \mathrm{u}^{(\ge k)} = \mathrm{f}_{\ge k}$ with
\begin{equation*}
\begin{split}
A_{\ge k}(\overline{\alpha_{k-1} \mathbf{j}_{\ge k}}, \overline{\alpha_{k-1}' \mathbf{j}_{\ge k}'}) & = \hat A_{<k}(\alpha_{k-1},\alpha_{k-1}') \A^{(k)}(j_k,j_k') \cdots \A^{(d)}(j_d,j_d'), \\
\mathrm{f}_{\ge k}(\overline{\alpha_{k-1} \mathbf{j}_{\ge k}}) & = \hat F_{<k}(\alpha_{k-1}) \f^{(k)}(j_k) \cdots \f^{(d)}(j_d),
\end{split}
\end{equation*}
where only the coordinates $<k$ are reduced.
Then, before Line \ref{al:ac:qr} in Alg.~\ref{alg:als-cross},
the \emph{residual} of the semi-reduced system is approximated in TT format with small TT ranks,
$$
\mathrm{z}_{\ge k} = \mathrm{f}{\ge k} - A_{\ge k} \mathrm{u}^{(\ge k)}, \quad \mathrm{z}_{\ge k}(\overline{\alpha_{k-1} \mathbf{j}_{\ge k}}) \approx Z^{| k \rangle}(\overline{\alpha_{k-1} j_k}) \mathbf{z}^{(k+1)}(j_{k+1}) \cdots \mathbf{z}^{(d)}(j_d),
$$
before \emph{expanding} the TT block $U^{|k\rangle}$ and orthogonalizing it to obtain $\left[U^{|k\rangle}
  \; Z^{|k\rangle}\right] = Q R_{\text{exp}}$.
To update the next TT block $U^{\langle k+1 |}$ consistently, we remove the last $\mathrm{rank}(Z^{|k\rangle})$ columns from $R_{\text{exp}}$ to obtain $R$ and then compute $R U^{\langle k+1 |}$.

To compute $Z^{| k \rangle}$ efficiently, we use the same hybrid ALS-cross procedure, but on different index sets, tailored to the residual.
Applying the recursive maxvol algorithm to the TT blocks of $\mathbf{z}^{(k)}$, we obtain an index set $\mathcal{\tilde J}_{>k}$, subsample the matrix, right hand side and the solution,
\begin{equation*}
\begin{split}
\tilde A_k & = \sum_{\gamma_{k-1},\gamma_k=1}^{R_{k-1},R_k}\hat A_{<k}(:,:,\gamma_{k-1}) \otimes \mathbf{A}^{(k)}_{\gamma_{k-1},\gamma_k} \otimes A^{(>k)}_{\gamma_k}(\mathcal{\tilde J}_{>k},\mathcal{\tilde J}_{>k}), \\
\mathrm{\tilde f}_k & = \sum_{\delta_{k-1},\delta_k=1}^{\rho_{k-1},\rho_k}\hat F_{<k}(:,\delta_{k-1}) \otimes \mathbf{f}^{(k)}_{\delta_{k-1},\delta_k} \otimes F^{(>k)}_{\gamma_k}(\mathcal{\tilde J}_{>k}),
\end{split}
\end{equation*}
$\tilde U^{|k\rangle}  = U^{|k\rangle} \big[U^{(>k)}(:, \mathcal{\tilde J}_{>k})\big],$ and update the TT block of the residual as $\mathrm{z}^{(k)} = \mathrm{\tilde f}_k - \tilde A_k \mathrm{\tilde u}^{(k)}$.

\begin{remark}
\label{rem:rank-increase}
The number of AMEn iterations to obtain convergence depends on two factors.
Firstly, if the initial TT ranks are underestimated, we need to raise them to achieve the desired accuracy.
The number of iterations in this stage is equal to the difference in the solution ranks, divided by the rank of $Z^{|k\rangle}$.
Secondly, once the proper ranks are achieved, a few more iterations
might be necessary for the solution to actually converge.
For the parametric PDE under consideration we have found that a single iteration in this second stage is sufficient.
Since in the log-normal case the ranks of the solution and of the
coefficient turn out to be approximately equal, this makes the treatment of the
log-normal case particularly simple. Using the coefficient as the
initial guess, Alg.~\ref{alg:als-cross} needs only one iteration overall.
In the affine case, the TT ranks of the solution can be several times larger than the ranks of the coefficient.
It is still advantageous though, to start with the coefficient and to conduct several iterations to increase the ranks.
\end{remark}

\subsection{Complexity of the algorithm}

Let us now consider the cost of each of the steps in Alg.~\ref{alg:als-cross}.
For brevity, assume that there exists an $r \in \mathbb{N}$, such that $R_k \le r$ and $r_k \le r$, for all $k=0,\ldots,d$.
Lines \ref{al:ac:qr} and \ref{al:ac:lq}--\ref{al:ac:Jr} of Alg.~\ref{alg:als-cross} (QR, maxvol decompositions and updates of index sets) require $\mathcal{O}(nr^3)$ operations.
Lines \ref{al:ac:evl} and \ref{al:ac:evr} (updates of TT blocks for $k>0$) require the construction and solution of $\le nr$ dense systems of size $\le r \times r$ at a cost of $\mathcal{O}(nr^4)$ each.
The computation of the projections in Line \ref{al:ac:Al} has the same complexity.
In particular, each of the matrix projections in \eqref{eq:Aleft} requires computing three tensor-tensor products.
The projection of the right-hand side \eqref{eq:flocal} has a cost of $\mathcal{O}(nr^3)$ which can be neglected.
Therefore, the TT blocks for $k=1,\ldots,d$ can be computed at the cost not greater than $\mathcal{O}(dnr^4)$ per iteration
with the hidden constant being of order 1.
This includes the products in the construction of \eqref{eq:PGmat}, the direct inversions of the dense blocks, and the products in the update of \eqref{eq:Aleft}.

We need to multiply a $N \times r$ by a $r \times r$ matrix to
construct the deterministic coefficients for~\eqref{eq:A_f_fix}, resulting in a cost of order
$\mathcal{O}(N r^2)$ for the setup of Line~\ref{al:ac:det}.
The cost of the solution of the deterministic problems in Line~\ref{al:ac:det} depends on the particular model and solver.
For an optimal solver (e.g. multigrid), this step has
$\mathcal{O}(Nr)$ complexity, in general it can be of $\mathcal{O}(K N^\lambda r)$ with $\lambda>1$ and/or $K\gg 1$.
The deterministic solutions in Line
\ref{al:ac:qr0} can be orthogonalized in $\mathcal{O}(Nr^2)$ operations, while the
deterministic matrices can be projected in $\mathcal{O}(Nr^3)$ operations.

The total cost per iteration (with the hidden constants being of order $1$) is therefore
\begin{equation}
\mathcal{O}(K N^\lambda r + Nr^3 + dnr^4).
\label{eq:cost}
\end{equation}
Typical magnitudes for the parameters in \eqref{eq:cost} are $N \sim 10^4-10^5$ (for a two-dimensional physical space), $d \sim 20-100$, $r \sim 100-500$, and $n \sim 5-10$ (for a Gaussian grid).
For smaller ranks, the complexity is mainly due to the deterministic part with the cost proportional to $N$.
However, for larger ranks the $r^4$-term becomes dominant, and the stochastic part can consume most of the time.
We will explore this in detail via some concrete numerical experiments.

\section{Numerical examples}
\label{sec:numer_main}
We consider the following instance of the PDE problem \eqref{eq:pdeproblem}: we choose $\Omega = [0,1]^2$, the right-hand side $f=0$, and the boundary conditions
\begin{equation}
u|_{x^1=0}=1, \quad u|_{x^1=1}=0, \quad \left.\frac{\partial u}{\partial n}\right|_{x^2=0} = \left.\frac{\partial u}{\partial n}\right|_{x^2=1} = 0
\label{eq:bc}
\end{equation}
This allows us to introduce, for any $m \in \mathbb{N}$, a uniform rectangular partitioning
$$
\mathcal{T}_h = \left\{[x^1_{i_1}, x^1_{i_1+1}] \times [x^2_{i_2}, x^2_{i_2+1}]: x^1_{i_1} =  i_1/m \mbox{ and } x^2_{i_2} = i_2/m, \quad i_1,i_2=0,\ldots,m-1\right\},
$$
with mesh width $h=1/m$ and the corresponding finite element space
$$
V_h = \text{span} \left\{\phi_i(x): i=\overline{i_1 i_2}=1,\ldots,N \right\}, \qquad N=\left(m-1\right)\left(m+1\right),
$$
where $i_1 \in \{1,\ldots,m-1\}$, $i_2 \in \{0,\ldots,m\}$ and $\phi_i(x)$ is a continuous, piecewise bilinear function on $[0,1]^2$, which is $1$ at $(i_1h, i_2h)$ and $0$ at all other vertices of the uniform grid.

The coefficient is based on the following expansion, proposed in \cite{Schneider-HT-SFEM-2016},
\begin{equation}
\begin{split}
w(x,\mathbf{y}) & = \sum_{k=1}^{d} y^k \psi_k(x), \quad \text{where}  \ \ \psi_k(x) = \sqrt{\eta_k} \cos(2\pi \rho_1(k) x^1) \cos(2\pi \rho_2(k) x^2), \\
\rho_1(k) & = k - \tau(k)\frac{(\tau(k)+1)}{2}, \quad \rho_2(k) = \tau(k)-\rho_1(k), \quad \tau(k)  = \left\lfloor -\frac{1}{2} + \sqrt{\frac{1}{4}+2k} \right\rfloor, \\
\eta_k & = \sigma^2 \frac{\mathcal{D}_k}{\sum_{m=1}^{d} \mathcal{D}_m}, \quad \text{and} \ \  \mathcal{D}_k = \left\{\begin{array}{ll}1, & k\le k_0, \\ (k-k_0)^{-\nu-1}, & k> k_0. \end{array}\right.
\end{split}
\label{eq:kle_art}
\end{equation}
Here, $d$ is chosen such that $\|\psi_{d+1}\|_{\infty}$ is smaller than the required tolerance.
This coefficient mimics the Karhunen-Lo\`eve expansion (KLE) induced by the Matern covariance function with smoothness parameter $\nu$, variance $\sigma^2$ and correlation length $\ell=1/k_0$.
Unless otherwise specified, we use $\sigma^2=\ell=1$ in the experiments.
The coefficient is constructed in either \emph{affine} or \emph{log-}form,
$$
c(x,\mathbf{y}) = 10+w(x,\mathbf{y}), \quad \mbox{or} \quad c(x,\mathbf{y}) = \exp\left(w(x,\mathbf{y})\right),
$$
respectively. To ensure well-posedness, only uniform $\mathbf{y}$ are used in the affine case,
while both uniformly and normally distributed parameters are used in the log-case.

The spatial grid levels vary from $1$ to $5$, with mesh size halved from level to level, starting from $h=1/32$ at Level $1$.
The stochastic parameters $y^k$ are discretized with $n$ Gauss-Hermite or Gauss-Legendre points for normal or uniform $\mathbf{y}$, resp.
By default, we use $n=7$ for the first dimension, and coarsen this
anisotropically in the latter dimensions (see Sec.~\ref{sec:ani} for details).

The (vector-valued) QoI consists of $10$ moments of the solution
averaged over a subdomain,
\begin{equation}
Q_p = \mathbb{E}\left[\mathbf{Q}^p(\mathbf{y})\right], \quad p=1,\ldots,10, \qquad \mathbf{Q}(\mathbf{y}) = \int_{[\frac{6}{8}, \frac{7}{8}]\times[\frac{7}{8}, \frac{8}{8}]} u(x,\mathbf{y}) dx - 0.2,
\label{eq:QoI}
\end{equation}
where the offset by $0.2$ is added to make the magnitudes of the first two moments comparable.
We show how these moments can be used for estimating the PDF of $\mathbf{Q}$ in Section \ref{sec:num_PDF}.

The TT solver is carried out as follows.
Since $10+w(x,\mathbf{y})$ is affine in $\mathbf{y}$, it is constructed straightforwardly in the TT format.
In the log-case, we approximate the exponential of $w$ at all spatial and parametric points using the TT-Cross algorithm (Alg.~\ref{alg:ttcross}).
After the full solution is computed via the ALS-Cross algorithm (Alg.~\ref{alg:als-cross}), integration over $x$ in \eqref{eq:QoI} involves only the zeroth (spatial) TT block.
The resulting quantity depends only on $\mathbf{y}$, and is raised to the $p$-th power via the TT-Cross algorithm. The moments of $\mathbf{Q}$ are computed by multiplying the obtained tensor with the quadrature weights (the latter admit a rank-1 TT decomposition).
To compute the coefficient and the QoI, we initialize Alg.~\ref{alg:ttcross} with $800$ random indices $\mathcal{J}_{>k}$, for each $k=0,\ldots,d-1$.
This is a reasonable trade-off between overhead and reliability, since optimal TT ranks range from $50$ to $500$, depending on the model parameters.

For the ALS-Cross algorithm we have two strategies, depending on the form of the coefficient, see Rem.~\ref{rem:rank-increase}.
In the log-normal and log-uniform case, the ranks of the coefficient and the solution are comparable,
and we can use the coefficient as an initial guess.
In the affine case, we can start from the coefficient and increase the ranks in the course of the iteration, or we can adopt the strategy used for the coefficient, i.e., to start from a high-rank initial guess and to truncate the ranks to the optimal values.
However, the rank of the initial guess governs the number of deterministic problems we need to solve;
if it is too large, the overall complexity will be too high.
We compare these two strategies below.

Let us now briefly describe the methods we benchmark our ALS-Cross algorithm against.
The first attempts of applying tensor techniques to stochastic PDEs used classical Richardson iteration, preconditioned by the mean-field operator \cite{khos-pde-2010,KhSch-Galerkin-SPDE-2011}.
A slightly improved method is the preconditioned steepest descent iteration where linear combinations of vectors, scalar and matrix products are replaced by their TT counterparts, followed by rank truncation. We will benchmark against this method in the affine case. The mean-field preconditioner is efficient for small
variances of the random field, but in general, excessive TT ranks of the preconditioned residual introduce a significant overhead. One way to reduce the complexity of classical iterative methods is to avoid multiplying the TT approximations of $A$ and $u$ and to compute directly the approximate residual using adaptive truncation both in the TT ranks and in the parametric resolution \cite{bachmayr-sparse-or-lr-2017},
but we will not compare against that approach.

Quasi Monte Carlo (QMC) methods are equal-weight quadrature rules
where the nodes are chosen deterministically. Here, we use a randomised
deterministic lattice rule, constructed by minimising the worst-case error
component by component \cite{graham-QMC-2011}.
We can also leverage the hierarchy of spatial
discretizations for the PDE problem and design a multilevel QMC (MLQMC) algorithm
\cite{Scheichl-mlqmc-lognorm-2017}.
In MLQMC, large numbers of samples are only needed on
the coarsest spatial grid.
The solutions on finer grids are then used to correct this coarse-grid solution.
When the difference between the solutions on two consecutive grids is
small, only a few QMC samples of the fine-grid solution are
sufficient to deliver the same accuracy as the single-level approach.

Finally, the Sparse Grid approach \cite{griebel-sparsegrids-2004} is a particular realization of stochastic collocation,
where instead of the full tensor product grid
we use a union of grid surpluses with bounded total polynomial degree.
Individual degrees for each variable can be weighted (anisotropically) to gain additional sparsity.
For the stochastic PDE, one can introduce lower degrees in the latter,
less important dimensions and the weights can be determined adaptively \cite{Gerstner-adapt-SG-2003}.

The TT algorithms\footnote{The codes are available from \url{http://people.bath.ac.uk/sd901/als-cross-algorithm}.} were implemented within the TT-Toolbox \cite{tt-toolbox}.
The QMC points are produced via randomly shifted lattice rules \cite{Scheichl-mlqmc-lognorm-2017}, using the generating vector from \cite[\texttt{lattice-39102-1024-1048576.3600.txt}]{Kuo-lattice-url}.
Sparse grid experiments are carried out using an adaptive sparse grid algorithm from the \texttt{spinterp} toolbox \cite{spinterp,spinterp-doc}.
This algorithm is only applicable for uniform $\mathbf{y}$ though, since it employs a Clenshaw-Curtis grid with boundary points, giving zero or infinite coefficients for normally distributed $\mathbf{y}$.
 The CPU times, measured in seconds, are summed over all runs.
The numerical results were obtained in Matlab R2015b on one node of
the University of Bath \texttt{Balena} cluster with a 2.6GHz Intel Xeon processor.


\subsection{Equilibration of different errors in the approximation}
\label{sec:ani}

We are going to benchmark the algorithms with respect to the total modelling error. To optimize the performance and to guarantee a fair comparison, we aim to equilibrate all the errors introduced at different stages: truncation of the coefficient expansion, discretization in space and in the parameters, approximation in TT format and quadrature errors in the QMC method. In this section, we investigate some of these individual errors, study their dependence on key modelling parameters and describe how we equilibrate them.

We begin with the truncation error.
An exact parametrization of a continuous random field $c(x,\omega)$, e.g., a Karhunen-Lo\`eve expansion of a Gaussian random field with Mat\`ern covariance, involves an infinite number of parameters, $c(x,y^1,y^2,\ldots)$.
For practical purposes, we have to truncate such an expansion to a finite number of parameters, and thus introduce an error.
In the case of the affine expansion $c(x,\mathbf{y}) = 10+w(x,\mathbf{y})$ with $w$ defined in \eqref{eq:kle_art}, the norm of $\psi_k$ is bounded by the function $\mathcal{D}_k$, which is asymptotically proportional to $k^{-\nu-1}$.
This decaying sequence can be used for estimating the truncation dimension and the corresponding error. A similar dependence of the truncation error on the norm of the functions $\psi_k$ and thus on $k^{-\nu-1}$ can also be shown in the log-normal case \cite{charrier-trunc-lognorm-2013}. In both cases, we truncate the sum by ensuring that the first discarded term has a norm less than a given tolerance $\delta$, i.e. $\|\psi_{d+1}\|_{\infty} \le \delta$.

Next, let us study the error due to the stochastic collocation method and its dependence on the maximal polynomial degree $n$ (and thus the number of collocation points) in any parameter direction. We consider both a uniform grid, where the same $n_k=n$ is used for all components $y^k$, $k=1,\ldots,d$, as well as a heuristic anisotropic coarsening strategy, reducing the polynomial degree as $k$ increases. The decay of $\mathcal{D}_k$ and thus of the norm of the $k$th term in the expansion of $w$ suggests the use of fewer discretization points for the latter $y^k$, i.e. an anisotropic grid. In the TT format, it is possible to estimate the quadrature error in computing the QoI for each individual term and we can bound the total error due to the choice of polynomial degree in the stochastic collocation approximation by
\begin{equation}
\widehat{\eps}_{\text{SC}} \le C \sum_{k=1}^{d} \mathcal{D}_k \widehat{\eps}_q(n_k),
\label{eq:err_ani}
\end{equation}
where $\widehat \eps_q(n_k)$ is the quadrature error in the $k$-th variable.
Our anisotropic coarsening strategy aims to balance the terms in
\eqref{eq:err_ani} such that
$\mathcal{D}_k \widehat \eps_q(n_k) = \delta/d$, where $\delta$ is again the prescribed error tolerance. For analytic functions, Gauss quadrature converges exponentially,
i.e., $\widehat \eps_q(n_k) = \exp(-cn_k)$, suggesting the choice $n_k = \frac{1}{c}\left(\log \mathcal{D}_k - \log (\delta/d)\right)$.
Denoting $n_1=n$ and setting $n_d=1$ for the last variable, we obtain
\begin{equation}
n_k = \left\lceil n + (1-n)\frac{\log\mathcal{D}_k}{\log \mathcal{D}_d} \right\rceil,
\label{eq:n_ani}
\end{equation}
where $n$ is chosen such that the total quadrature error is less than $\delta$.


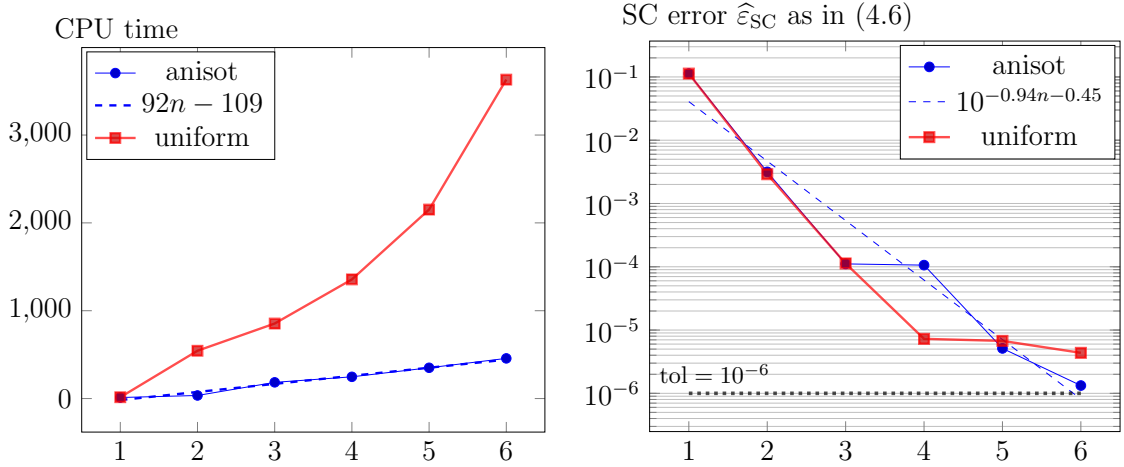
\begin{figure}[t]
\centering
\caption{CPU times (left) and errors (right) for different parameter grid sizes $n$, with the same number of points for different dimensions (uniform) and anisotropic coarsening (anisot)}
\label{fig:ny}
\resizebox{\figscale\linewidth}{!}{
\begin{tikzpicture}
  \begin{axis}[%
  xmode=normal,
  ymode=normal,
  ylabel=CPU time,
  legend style={at={(0.01,0.99)},anchor=north west},
  ]
  \addplot+[] coordinates{(1, 11.7086) (2, 34.8527) (3, 184.4629) (4, 247.3686) (5, 350.3241) (6, 457.2798)};\addlegendentry{anisot}; 
  \addplot+[blue,dashed,no marks,line width=1pt,domain=1:6] {92*x-109};\addlegendentry{$92 n-109$};
  \pgfplotsset{cycle list shift=-1};
  \addplot+[line width=1pt,opacity=0.7] coordinates{(1, 15.0161) (2, 544.0828) (3, 854.1654) (4, 1.3563e+03) (5, 2.1518e+03) (6, 3.6311e+03)};\addlegendentry{uniform}; 
  \end{axis}
 \end{tikzpicture}
}
\resizebox{\figscale\linewidth}{!}{
\begin{tikzpicture}
  \begin{axis}[%
  xmode=normal,
  ymode=log,
  ylabel=SC error $\widehat \eps_{\text{SC}}$ as in \eqref{eq:err_det},
  legend style={at={(0.99,0.99)},anchor=north east},
  ]
  \addplot+[] coordinates{(1, 1.1291e-01) (2, 3.1614e-03) (3, 1.1099e-04) (4, 1.0615e-04) (5, 5.1279e-06) (6, 1.3277e-06)};\addlegendentry{anisot}; 
  \addplot+[blue,dashed,no marks,domain=1:6] {10^(-0.94*x-0.45)};\addlegendentry{$10^{-0.94 n-0.45}$};
  \pgfplotsset{cycle list shift=-1};
  \addplot+[line width=1pt,opacity=0.7] coordinates{(1, 0.1129) (2, 0.0029) (3, 1.1282e-04) (4, 7.2587e-06) (5, 6.7424e-06) (6, 4.3642e-06)};\addlegendentry{uniform}; 
  \addplot+[black,opacity=0.7,line width=1.5pt,no marks,dotted,domain=1:6] {1e-6}; \node[anchor=south west] at (axis cs:0.5,1.0e-6) {\footnotesize $\mathrm{tol}=10^{-6}$};
  \end{axis}
 \end{tikzpicture}
}
\end{figure}

In Fig.~\ref{fig:ny}, we investigate numerically the quadrature error due to collocation in the stochastic parameters $\mathbf{y}$.
We fix the spatial mesh level at $4$ (i.e. $h=1/256$), let $\nu=3$ and set the TT approximation and KLE truncation thresholds to $\delta = 10^{-6}$, varying the numbers of grid points $n_1,\ldots,n_d$ in $\mathcal{Y}_n$.
We consider both the uniform grid with $n_k=n$, for all $k=1,\ldots,d$, as well as the anisotropic coarsening strategy in \eqref{eq:n_ani}, discussed above.
The total stochastic collocation error $\widehat \eps_{\text{SC}}$ is estimated by considering the relative error in the Frobenius-norm of the vector of moments \eqref{eq:QoI}:
\begin{equation}
\widehat \eps_{\text{SC}} = \|Q - Q_\star\|_2/\|Q_\star\|_2\,
\label{eq:err_det}
\end{equation}
where the reference solution $Q_{\star}$ is computing using all the same model parameters except for a uniform collocation grid with $n=11$.
We see that on average, the errors for the uniform and for the anisotropic grid are the same, while the CPU time for the uniform grid is significantly larger, and moreover, grows faster than linearly due to the growth of the TT ranks.
In the anisotropic setup, the growth is perfectly linear in $n$.

\begin{figure}[t]
\centering
\caption{Spatial discretization error (left) and QMC error (right)
for log-normal field \eqref{eq:kle}.}
\label{fig:eps-h}
\resizebox{\figscale\linewidth}{!}{
 \begin{tikzpicture}
  \begin{axis}[%
  ymode=log,
  xmode=normal,
  xtick={1,2,3,4},
  xlabel=$\ell \equiv \log_2 (1/h) - 4$,
  ylabel=FE error $\widehat \eps_{\text{FE}}$ similar to \eqref{eq:err_det},
  ]
  \addplot+[] coordinates{(1,     1.2272e-03) (2,     3.0344e-04) (3,     7.1396e-05) (4,     1.7662e-05)};\addlegendentry{$\nu=3$};
  \addplot+[] coordinates{(1,     1.9312e-03) (2,     4.6331e-04) (3,     9.1921e-05)};\addlegendentry{$\nu=1.5$};
  \addplot+[no marks,domain=1:4] {2^(-2.034*x-7.613)};\addlegendentry{$2^{-2.034 \cdot \ell-7.613}$};
  \end{axis}
 \end{tikzpicture}
}
\resizebox{\figscale\linewidth}{!}{
 \begin{tikzpicture}
  \begin{axis}[%
  ymode=log,
  xmode=normal,
  xlabel=$\ell \equiv \log_2 \mathcal{N}$,
  ylabel=QMC error $\widehat \eps_{\text{QMC}}$ similar to \eqref{eq:err_det},
  ]
  \addplot+[] coordinates{(8   , 2.6836e-02) (9   , 1.4105e-02) (10  , 8.0864e-03) (11  , 5.8352e-03) (12  , 2.5918e-03) (13  , 1.0165e-03) (14  , 6.3514e-04)};\addlegendentry{$\nu=3$};

  \addplot+[] coordinates{(8, 2.8041e-02) (9, 2.1009e-02) (10, 7.5942e-03) (11, 3.9207e-03) (12, 1.4598e-03) (13, 1.0636e-03)};\addlegendentry{$\nu=1.5$};
  \addplot+[no marks,domain=8:14] {2^(-0.971*x+2.885)};\addlegendentry{$2^{-0.971 \cdot\ell+2.885}$};
  \end{axis}
 \end{tikzpicture}
}
\end{figure}
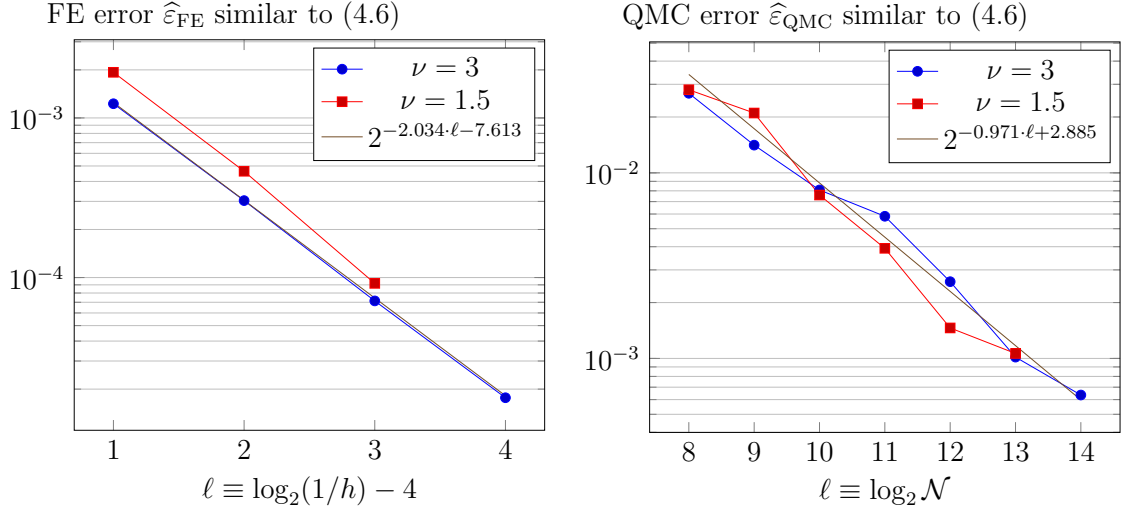

As stated above, for $\nu \ge 1$, the spatial discretization error has an algebraic
decay rate of $\mathcal{O}(h^2)$.
We estimate the precise dependence numerically for different values of $\nu$ by comparing again the vector of moments with the reference solution as in \eqref{eq:err_det}, see Fig.~\ref{fig:eps-h} (left).
The moments are computed using the same $\mathcal{N}=2^{20}$ QMC points for all levels and a tolerance of $\delta=10^{-5}$ for the truncation error,
with the reference solution $Q_\star$ computed with a spatial mesh size of $h_* = 1/1024$.
In the same way, the QMC quadrature error for different numbers of samples (with all other parameters fixed) is shown in Fig.~\ref{fig:eps-h} (right).
Here, the spatial mesh size is fixed to $1/256$, and the reference solution is computed with $\mathcal{N}=2^{20}$ QMC points.
We see that the spatial and QMC errors are not sensitive to the choice of $\nu$.



We finish with an analysis of the TT approximation errors.
The TT algorithms are parametrized again by the same stopping tolerance $\delta$; the TT ranks are adapted such that the relative mean error in the Frobenius vector norm is below that tolerance.
In Tab.~\ref{tab:tterr}, we vary $\delta$ from $10^{-2}$ to $10^{-5}$ and compute Monte Carlo (MC) estimates of the expected values of the relative errors in the coefficient $c$ and  the solution $u$ of \eqref{eq:pdeproblem}. In particular, we compute
\begin{equation}
\widehat \eps_{\text{TT,c}} = \frac{1}{\mathcal{N}} \sum_{i=1}^{\mathcal{N}} \frac{\left\| \c_{\text{TT}}(:,\mathbf{y}_i) - \c(:,\mathbf{y}_i) \right\|_{\infty}}{\left\| \c(:,\mathbf{y}_i) \right\|_{\infty}} \ \ \text{and} \ \
 \widehat \eps_{\text{TT,u}} = \frac{1}{\mathcal{N}} \sum_{i=1}^{\mathcal{N}} \frac{\left\| \u_{\text{TT}}(:,\mathbf{y}_i) - \u(:,\mathbf{y}_i) \right\|_{2}}{\left\| \u(\cdot,\mathbf{y}_i) \right\|_{2}},
 \label{eq:errtt}
\end{equation}
where $\c_{\text{TT}}$ and $\u_{\text{TT}}$ are the TT approximations of the exact coefficient and the solution for the case of the log-normal coefficient field $c(x,\mathbf{y}_i) = \exp(w(x,\mathbf{y}_i))$ with $\nu=3$, and where $\{\mathbf{y}_i:i=1,\ldots,\mathcal{N}\}$ is a set of $\mathcal{N}$ i.i.d. standard normals. The norm used for $\widehat \eps_{\text{TT,c}}$ approximates the error in the $L^{\infty}$ norm, while the norm used for $\widehat \eps_{\text{TT,u}}$ approximates the error in the $L^2$ norm.
The total number of samples $\mathcal{N}=10^4$, the spatial mesh size is set to $h=1/256$, $n=7$ and the truncation tolerance is $10^{-5}$.
We see that the actual errors decay (roughly) linearly with $\delta$, and the $L^2$ norm of the error is within the requested tolerance.

\begin{table}[ht]
 \centering
 \caption{Relative errors in the $L^\infty$-norm of the coefficient and in the $L^2$-norm of the solution with respect to the TT approximation tolerance (MC estimate with 95\% confidence interval).}
 \label{tab:tterr}

 \vspace*{1ex}
 \begin{tabular}{c|cccc}
  TT tolerance  $\delta$  & $10^{-2}$  & $10^{-3}$  & $10^{-4}$  & $10^{-5}$   \\ \hline
  $\widehat \eps_{\text{TT,c}}$     & $\text{1.53e-2}\pm\text{3.3e-4}$ & $\text{2.15e-3}\pm\text{4.7e-5}$ & $\text{2.5e-4}\pm\text{5.7e-6}$ & $\text{2.92e-5}\pm\text{6.8e-7}$  \\
  $\widehat \eps_{\text{TT,u}}$     & $\text{4.91e-3}\pm\text{4.3e-5}$ & $\text{6.10e-4}\pm\text{1.1e-5}$ & $\text{8.7e-5}\pm\text{2.0e-6}$ & $\text{1.01e-5}\pm\text{1.9e-7}$  \\
 \end{tabular}
 \color{black}
\end{table}

\subsection{Benchmarking the ALS-Cross algorithm}

Let us now come to the main experiments, where we compare the TT scheme and the alternative techniques in terms of CPU time and number of deterministic solves versus total error. We vary the spatial grid level from 1 ($h=1/32$) to 5 ($h=1/512$). To equilibrate all errors, we use the estimate of the spatial discretization error in Fig.~\ref{fig:eps-h} for each grid level as the tolerance $\delta$ for the truncation error, for the stochastic collocation error and for the quadrature error (either due to QMC or TT approximation). This heuristic should ensure that the total error is less than $4\delta$.

Since both the TT and the QMC approaches involve randomness, we conduct $16$ runs of each and compare the results. The total error in the QoI is log-averaged over the runs,
\begin{equation}
\widehat \eps_{\text{tot}} = \exp\left(\frac{1}{16}\sum_{\iota=1}^{16} \log\|Q(\iota)-Q_\star\|_2 - \log\|Q_\star\|_2\right),
\label{eq:err}
\end{equation}
where $\iota$ enumerates the runs, and $Q_\star$ is the reference solution computed using MLQMC with $6$ levels, $h_L = 1/1024$ and all other errors equilibrated.
The $2$-norm of the vector of moments $Q$ is used since it governs the Kullback-Leibler divergence of the estimated PDF (see~\cite{Chernov-maxentr-2016}).

Since all algorithms are built upon state of the art Matlab software and run on the same computer, the CPU time should give a reliable complexity measure.
In addition, we also compare the methods in terms of number of deterministic PDE solves, since they constitute most of the cost, especially for complex PDEs. In fact, numerical experiments confirm that for our model problem, the solution of the $\mathcal{O}(r)$ deterministic PDEs (in Line \ref{al:ac:det} in Alg.~\ref{alg:als-cross}) dominates when the TT ranks remain relatively small (below $r=100$). This corresponds to the first term in the complexity estimate \eqref{eq:cost}. When $r$ is larger, the other steps in Alg.~\ref{alg:als-cross} dominate the cost, due to their $\mathcal{O}(r^3)$ and $\mathcal{O}(r^4)$ complexity. More specifically, in our tests with log-normal random fields on spatial grid level 4 ($h=1/256$), Line \ref{al:ac:det} in Alg.~\ref{alg:als-cross} constitutes 45\% of the cost for smoother fields ($\nu=3$) and about 5\% for less smooth fields ($\nu=1.5$).

\begin{remark}
\label{rem:norm_nsolve}
For MLQMC, the deterministic solves on different levels have different complexity, but we know (or can estimate) very accurately how this cost scales with grid level $\ell$. Thus, to allow a comparison with the other methods in terms of number of deterministic PDE solves, we
normalize the cost of each PDE solve on coarser levels, so that they represent the correct fraction of the cost of a PDE solve on the finest grid. The sum of these normalized costs over all the levels is then an accurate estimate of the equivalent cost in terms of number of PDE solves on the finest grid when compared to the other (single-level) methods.
\end{remark}

\subsubsection{Affine coefficients}

\begin{figure}[t]
\centering
\caption{CPU time (left) and number of deterministic PDE solves (right)
vs. $\log_{10}\widehat \eps_{\text{tot}}$ for the affine-uniform coefficient
\eqref{eq:kle_aff} with $\nu=3$ (\# solves in mlqmc computed according to Rem.~\ref{rem:norm_nsolve}).}
\label{fig:nu4_af}
\resizebox{\figscale\linewidth}{!}{
\begin{tikzpicture}
  \begin{axis}[%
  xmode=normal,
  ymode=log,
  ylabel=CPU time,
  legend style={at={(0.01,0.01)},anchor=south west},
  x filter/.code={\pgfmathparse{log10(\pgfmathresult*1)}\pgfmathresult},
  ]
   \addplot+[] coordinates{
                           (5.1461e-05, 8.2192e+01)
                           (1.1564e-05, 1.3668e+02)
                           (2.5281e-06, 5.9923e+02)
                           (6.1759e-07, 2.8371e+03)
                           }; \addlegendentry{tt1}
   \node[anchor=west] at (axis cs:-6.0,2e3) {\footnotesize\color{blue} 1.03};

   \addplot+[blue,dashed,mark options={blue},mark=otimes,line width=1pt]  coordinates  {
        (3.4181e-05, 4.5817e+01)
        (1.1170e-05, 2.3575e+01)
        (3.0299e-06, 1.5396e+02)
        (7.3402e-07, 9.6605e+02)
        (1.1733e-07, 1.0052e+04)
   };\addlegendentry{ttK}
   \node[anchor=east] at (axis cs:-6.5,2e3) {\footnotesize\color{blue} 1.32};

   \addplot+[green!50!black,mark options={green!50!black},mark=triangle*] coordinates  {
            (2.3287e-05, 1.1434e+02)
            (7.4780e-06, 3.5380e+02)
            (3.4529e-06, 3.2237e+03)
            (3.0144e-06, 1.9785e+04)
   };\addlegendentry{psd}

   \pgfplotsset{cycle list shift=-1}
   \pgfplotsset{cycle list shift=-2}

   \addplot+[]  coordinates{
       (2.8927e-05, 5.1051e+02)
       (7.2455e-06, 9.1791e+03)
       (2.0296e-06, 1.6696e+05)
   };\addlegendentry{qmc}
   \node[anchor=west] at (axis cs:-5.6,1e5) {\footnotesize\color{red} 2.18};

   \pgfplotsset{cycle list shift=-1}
   \addplot+[]          coordinates  {
        (2.1249e-05, 2.0835e+03)
        (6.0704e-06, 1.3445e+04)
        (1.8531e-06, 5.6245e+04)
        (3.6645e-07, 2.5615e+05)
   };\addlegendentry{mlqmc}
   \node[anchor=west] at (axis cs:-6.3,2e5) {\footnotesize\color{black} 1.18};
  \end{axis}
 \end{tikzpicture}
}
\resizebox{\figscale\linewidth}{!}{
 \begin{tikzpicture}
  \begin{axis}[%
  xmode=normal,
  ymode=log,
  ylabel=\# solves,
  x filter/.code={\pgfmathparse{log10(\pgfmathresult*1)}\pgfmathresult},
  ]
   \addplot+[]  coordinates{
        (5.1461e-05, 6400)
        (1.1564e-05, 6400)
        (2.5281e-06, 6400)
        (6.1759e-07, 6400)
   }; 

   \addplot+[blue,dashed,mark options={blue},mark=otimes,line width=1pt]  coordinates{
          (3.4181e-05,   382)
          (1.1170e-05,   553)
          (3.0299e-06,  1287)
          (7.3402e-07,  1794)
          (1.1733e-07,  3603)
   }; 

   \pgfplotsset{cycle list shift=-1}
   \addplot+[]  coordinates{
        (2.8927e-05,  131072)
        (7.2455e-06,  524288)
        (2.0296e-06, 2097152)
   }; 

   \pgfplotsset{cycle list shift=0}
   \addplot+[]  coordinates{
        (2.1249e-05,   524224  )
        (6.0704e-06, 6.0405e+05)
        (1.8531e-06, 5.7174e+05)
        (3.6645e-07, 5.2631e+05)
   }; 
  \end{axis}
 \end{tikzpicture}
}
\end{figure}

We start with the affine coefficient and a fairly smooth random field ($\nu=3$).
In Fig.~\ref{fig:nu4_af}, we compare our new ALS-Cross algorithm with a simple preconditioned steepest descent iteration in the TT format (``{psd}''), as well as with the QMC methods, in single (``qmc'') and multilevel (``{mlqmc}'') form \cite{Scheichl-mlqmc-lognorm-2017}.
The ALS-Cross algorithm is invoked in two regimes: either one iteration, starting from a random initial guess with TT ranks $800$ (``{tt1}''), or $K$ iterations, starting from the coefficient as the initial guess (``{ttK}'').
In the latter case, the algorithm is stopped when the relative difference between the solutions at iterations $K$ and $K-1$ is less than the chosen TT tolerance $\delta$.
The numbers in the left figure denote slopes of least squares linear fits (discarding the right-most point(s) where necessary).

We see that the ALS-Cross schemes are clearly faster and more accurate than the other techniques.
The multilevel QMC method exhibits the same asymptotic slope, but its absolute run time is two orders of magnitudes larger.
Note that the minimal possible slope of cost with respect to error is $1$ for the deterministic problem, since the FE error decays as $\mathcal{O}(h^2)$, while the cost grows at least with $\mathcal{O}(h^{-2})$.
In the ALS-Cross and in the multilevel QMC methods this is indeed almost achieved for the overall cost.
The single-level QMC exhibits a slope of~$2$, as was also observed previously \cite{Scheichl-mlqmc-lognorm-2017}.
The simple preconditioned steepest descent method becomes too slow due to the large TT ranks of the residual at the latter iterations.
Its cost seems to be growing faster than algebraically.
For the affine case, the adaptive Richardson iteration in \cite{bachmayr-sparse-or-lr-2017} may allow to overcome this problem, but for larger variances such as those in Sec.~\ref{sec:lognormal} we still expect the number of iterations to be too high to be really competitive.

Comparing the two versions of the ALS-Cross method (``tt1'' and ``ttK''), we confirm that for lower accuracies, it is faster to perform several iterations, and to approach the TT ranks of the solution from below, as this requires fewer deterministic solves, which dominate the total cost.
However, the method works also very well when choosing a high initial rank and performing only one iteration of the ALS-Cross algorithm. In our experiments, the complexity of ``tt1'' approaches the complexity of ``ttK'' for higher values of the TT tolerance $\delta$, since in that case the difference between initial and final TT ranks is smaller.


\subsubsection{Log-normal and log-uniform coefficients}
\label{sec:lognormal}

\begin{figure}[t]
\centering
\caption{CPU time (left) and number of deterministic solves (right)
  vs.~$\log_{10}\widehat \eps_{\text{tot}}$ for log-uniform (top) and log-normal (bottom)
  coefficients with $\nu=3$ (\# solves in mlqmc computed according to Rem.~\ref{rem:norm_nsolve}).}
\label{fig:nu4}
\resizebox{\figscale\linewidth}{!}{
\begin{tikzpicture}
  \begin{axis}[%
  xmode=normal,
  ymode=log,
  ylabel=CPU time,
  legend style={at={(0.01,0.01)},anchor=south west},
  x filter/.code={\pgfmathparse{log10(\pgfmathresult*1)}\pgfmathresult},
  ]
   \addplot+[] coordinates{
            (3.3427e-03,   2.7478e+01)
            (1.3886e-03,   3.1553e+01)
            (2.7797e-04,   1.4613e+02)
            (3.9122e-05,   8.2544e+02)
            (7.2564e-06,   5.6563e+03)
   };\addlegendentry{tt}
   \node[anchor=east] at (axis cs:-4.6+0.6,3e2) {\footnotesize\color{blue} 0.98};

   \addplot+[]  coordinates{
        (1.3266e-03,   5.1263e+02)
        (3.8172e-04,   9.2032e+03)
        (8.1076e-05,   1.6829e+05)
   };\addlegendentry{qmc}
   \node[anchor=west] at (axis cs:-3.9+0.6,5e3) {\footnotesize\color{red} 2.06};

   \addplot+[] coordinates{
        (1.4151e-03,   1.6393e+01)
        (4.8741e-04,   2.3852e+02)
        (2.1791e-04,   6.1623e+03)
   };\addlegendentry{sg}
   \node[anchor=east] at (axis cs:-3.5,1.3e3) {\footnotesize\color{brown} 3.13};

   \addplot+[]  coordinates{
           (1.4220e-03,   5.2973e+02)
           (3.5324e-04,   3.5052e+03)
           (8.4416e-05,   1.4574e+04)
           (1.9382e-05,   6.0073e+04)
           (3.1710e-06,   1.8826e+05)
   };\addlegendentry{mlqmc}
   \node[anchor=west] at (axis cs:-5.1+0.6,4e4) {\footnotesize\color{black} 1.10};
  \end{axis}
 \end{tikzpicture}
}
\resizebox{\figscale\linewidth}{!}{
 \begin{tikzpicture}
  \begin{axis}[%
  xmode=normal,
  ymode=log,
  ylabel=\# solves,
  x filter/.code={\pgfmathparse{log10(\pgfmathresult*1)}\pgfmathresult},
  ]

   \addplot+[]  coordinates{
            (3.3427e-03,    344)
            (1.3886e-03,    509)
            (2.7797e-04,    728)
            (3.9122e-05,   1005)
            (7.2564e-06,   1363)
   };

   \addplot+[]  coordinates{
        (1.3266e-03,   1.3107e+05)
        (3.8172e-04,   5.2429e+05)
        (8.1076e-05,   2.0972e+06)
   };

   \addplot+[] coordinates{
        (1.4151e-03,   3.2370e+03)
        (4.8741e-04,   1.4119e+04)
        (2.1791e-04,   8.1695e+04)
   };

   \addplot+[]  coordinates{
           (1.4220e-03,   1.3101e+05)
           (3.5324e-04,   1.5059e+05)
           (8.4416e-05,   1.2986e+05)
           (1.9382e-05,   1.3816e+05)
           (3.1710e-06,   1.1214e+05)
   };

%
  \end{axis}
 \end{tikzpicture}
}\\
\resizebox{\figscale\linewidth}{!}{
\begin{tikzpicture}
  \begin{axis}[%
  xmode=normal,
  ymode=log,
  ylabel=CPU time,
  legend style={at={(0.99,0.99)},anchor=north east},
  x filter/.code={\pgfmathparse{log10(\pgfmathresult*1)}\pgfmathresult},
  ]
   \addplot+[] coordinates{
        (5.5508e-03,   3.3727e+01)
        (8.8048e-04,   3.8130e+01)
        (2.7934e-04,   1.7484e+02)
        (5.0192e-05,   9.9976e+02)
        (5.9748e-06,   6.8684e+03)
   };\addlegendentry{tt}
   \node at (axis cs:-4.6+0.6,3e2) {\footnotesize\color{blue} 1.03};

   \addplot+[]  coordinates{
         (1.4190e-03,   5.1391e+02)
         (4.2696e-04,   9.1695e+03)
         (9.1712e-05,   1.6858e+05)
   };\addlegendentry{qmc}
   \node at (axis cs:-4.1+0.6,5e4) {\footnotesize\color{red} 2.11};

   \pgfplotsset{cycle list shift=1}

   \addplot+[]  coordinates{
        (1.3707e-03,   1.0529e+03)
        (3.7348e-04,   9.7985e+03)
        (9.5152e-05,   2.8362e+04)
        (1.6594e-05,   2.3967e+05)
   };\addlegendentry{mlqmc}
   \node at (axis cs:-5.1+0.6,2e5) {\footnotesize\color{black} 1.18};

  \end{axis}
 \end{tikzpicture}
}
\resizebox{\figscale\linewidth}{!}{
 \begin{tikzpicture}
  \begin{axis}[%
  xmode=normal,
  ymode=log,
  ylabel=\# solves,
  x filter/.code={\pgfmathparse{log10(\pgfmathresult*1)}\pgfmathresult},
  ]
   \addplot+[] coordinates{
        (5.5508e-03,   5.3000e+02)
        (8.8048e-04,   6.7400e+02)
        (2.7934e-04,   8.0800e+02)
        (5.0192e-05,   1.0860e+03)
        (5.9748e-06,   1.3980e+03)
   };

   \addplot+[]  coordinates{
         (1.4190e-03,   1.3107e+05)
         (4.2696e-04,   5.2429e+05)
         (9.1712e-05,   2.0972e+06)
   };

   \pgfplotsset{cycle list shift=1}

   \addplot+[]  coordinates{
        (1.3707e-03,   2.6208e+05)
        (3.7348e-04,   3.4725e+05)
        (9.5152e-05,   3.4278e+05)
        (1.6594e-05,   4.4552e+05)
   };

  \end{axis}
 \end{tikzpicture}
}
\end{figure}

In the next experiments (Fig.~\ref{fig:nu4}--\ref{fig:nu2.5}), we focus on the log-uniform and the log-normal coefficients. We carry out systematic comparisons of the (tt1 version of the) ALS-Cross method, starting with the coefficient as the initial guess (``tt''),  with the single and multilevel QMC methods (``qmc'', ``mlqmc''), as well as with the dimension adaptive sparse grid method (``sg'', for uniform $\mathbf{y}$ only).
As before, we compare total CPU times (left) and numbers of deterministic solves (right), plotted against the estimate of the total error $\widehat \eps_{\text{tot}}$ in \eqref{eq:err}.

For the smooth field ($\nu=3$, Fig.~\ref{fig:nu4}) the ALS-Cross method is the fastest approach, due to small TT ranks in both the spatial and stochastic TT blocks.
For spatial grid level $5$ ($h=1/512$) with normally distributed $\mathbf{y}$, for example, it is almost $100$ times faster than the nearest competitor, multilevel QMC.
In the log-uniform case, for lower accuracies, the sparse grid algorithm is the fastest method, since a low polynomial order (and hence a small number of samples) suffices.
However, due to the simultaneous refinement of the spatial and parametric grids (the latter also grows with the total dimension) the cost grows quite steeply.
This makes the sparse grid approach significantly slower for higher accuracies in this example.

The difference in numbers of solves is more dramatic.
The TT-based scheme needs to solve only $r_0$ deterministic problems, where
$r_0$ grows logarithmically with the tolerance $\delta$, and is almost independent of the spatial grid size. However, since the stationary 2D diffusion equation on the unit square is
relatively easy to solve, this is not directly transferable to the comparison in terms of cost versus error. The computation of the TT approximation of the stochastic part contributes significantly to the overall cost. For more difficult PDE models (e.g. the
time-dependent Navier-Stokes equations in three space dimensions) we expect the deterministic part to become the most time-consuming part.


\begin{figure}[t]
\centering
\caption{CPU time (left) and number of deterministic solves (right)
  vs.~$\log_{10}\widehat \eps_{\text{tot}}$ for log-uniform (top) and log-normal (bottom)
  coefficients with $\nu=1.5$  (\#  solves in mlqmc computed according to Rem.~\ref{rem:norm_nsolve}).}
\label{fig:nu2.5}
\resizebox{\figscale\linewidth}{!}{
\begin{tikzpicture}
  \begin{axis}[%
  xmode=normal,
  ymode=log,
  ylabel=CPU time,
  legend style={at={(0.01,0.01)},anchor=south west},
  x filter/.code={\pgfmathparse{log10(\pgfmathresult*1)}\pgfmathresult},
  ]
   \addplot+[]  coordinates{
       (5.7744e-03,   4.7973e+01)
       (1.5144e-03,   1.5113e+02)
       (3.0322e-04,   2.5423e+03)
       (7.2848e-05,   5.0923e+04)
   };\addlegendentry{tt}
   \node[anchor=east] at (axis cs:-3.6+0.6,2e2) {\footnotesize\color{blue} 1.61};

   \addplot+[]  coordinates{
        (3.1154e-03,   5.1498e+02)
        (8.1736e-04,   9.1276e+03)
        (1.8248e-04,   1.6464e+05)
   };\addlegendentry{qmc}
   \node[anchor=west] at (axis cs:-4.3+0.6,1e5) {\footnotesize\color{red} 2.00};

   \addplot+[] coordinates{
        (3.8782e-03,   1.6871e+02)
        (1.6633e-03,   3.7230e+03)
        (1.1474e-03,   1.3858e+04)
   };\addlegendentry{sg}
   \node[anchor=west] at (axis cs:-2.9,1e4) {\footnotesize\color{brown} 3.52};

   \addplot+[]  coordinates{
         (3.0411e-03,   1.1165e+03)
         (7.3520e-04,   6.9480e+03)
         (1.8172e-04,   3.2889e+04)
         (3.9452e-05,   2.2596e+05)
   };\addlegendentry{mlqmc}
   \node[anchor=west] at (axis cs:-5.0+0.6,3e5) {\footnotesize\color{black} 1.21};
  \end{axis}
 \end{tikzpicture}
}
\resizebox{\figscale\linewidth}{!}{
 \begin{tikzpicture}
  \begin{axis}[%
  xmode=normal,
  ymode=log,
  ylabel=\# solves,
  x filter/.code={\pgfmathparse{log10(\pgfmathresult*1)}\pgfmathresult},
  ]
   \addplot+[]  coordinates{
       (5.7744e-03,   7.6600e+02)
       (1.5144e-03,   1.3860e+03)
       (3.0322e-04,   2.6190e+03)
       (7.2848e-05,   4.6900e+03)
   };

   \addplot+[]  coordinates{
        (3.1154e-03,   1.3107e+05)
        (8.1736e-04,   5.2429e+05)
        (1.8248e-04,   2.0972e+06)
   };

   \addplot+[] coordinates{
        (3.8782e-03,   2.8673e+04)
        (1.6633e-03,   2.0004e+05)
        (1.1474e-03,   2.0000e+05)
   };

   \addplot+[]  coordinates{
         (3.0411e-03,   2.6208e+05)
         (7.3520e-04,   3.2654e+05)
         (1.8172e-04,   3.0487e+05)
         (3.9452e-05,   4.3153e+05)
   };

  \end{axis}
 \end{tikzpicture}
}\\
\resizebox{\figscale\linewidth}{!}{
\begin{tikzpicture}
  \begin{axis}[%
  xmode=normal,
  ymode=log,
  ylabel=CPU time,
  legend style={at={(0.01,0.01)},anchor=south west},
  x filter/.code={\pgfmathparse{log10(\pgfmathresult*1)}\pgfmathresult},
  ]
   \addplot+[]  coordinates{
      (8.0304e-03,   7.6430e+01)
      (1.3996e-03,   4.6438e+02)
      (3.7910e-04,   7.9516e+03)
      (8.4656e-05,   1.0856e+05)
   };\addlegendentry{tt}
   \node[anchor=west] at (axis cs:-3.0+0.6,2e2) {\footnotesize\color{blue} 1.63};

   \addplot+[]  coordinates{
        (3.0703e-03,   5.1596e+02)
        (9.2732e-04,   9.1083e+03)
        (2.2281e-04,   1.6476e+05)
   };\addlegendentry{qmc}
   \node[anchor=west] at (axis cs:-4.1+0.6,1e5) {\footnotesize\color{red} 2.20};

   \pgfplotsset{cycle list shift=1}
   \addplot+[]  coordinates{
        (2.8050e-03,   2.2681e+03)
        (7.4504e-04,   1.3382e+04)
        (1.6528e-04,   7.9383e+04)
        (3.5624e-05,   3.2086e+05)
   };\addlegendentry{mlqmc}
   \node[anchor=west] at (axis cs:-3.3+0.6,4e3) {\footnotesize\color{black} 1.14};
  \end{axis}
 \end{tikzpicture}
}
\resizebox{\figscale\linewidth}{!}{
 \begin{tikzpicture}
  \begin{axis}[%
  xmode=normal,
  ymode=log,
  ylabel=\# solves,
  x filter/.code={\pgfmathparse{log10(\pgfmathresult*1)}\pgfmathresult},
  ]
   \addplot+[]  coordinates{
      (8.0304e-03,   1.3600e+03)
      (1.3996e-03,   2.0950e+03)
      (3.7910e-04,   3.6250e+03)
      (8.4656e-05,   5.7540e+03)
   };

   \addplot+[]  coordinates{
        (3.0703e-03,   1.3107e+05)
        (9.2732e-04,   5.2429e+05)
        (2.2281e-04,   2.0972e+06)
   };

   \pgfplotsset{cycle list shift=1}
   \addplot+[]  coordinates{
        (2.8050e-03,   5.2422e+05)
        (7.4504e-04,   6.6057e+05)
        (1.6528e-04,   7.5627e+05)
        (3.5624e-05,   6.4913e+05)
   };

  \end{axis}
 \end{tikzpicture}
}
\end{figure}

For a less smooth field ($\nu=1.5$, Fig.~\ref{fig:nu2.5}), the slope for the cost of
the ALS-Cross method increases, in particular becoming larger than the respective slope of multilevel QMC. However, in absolute terms the ALS-Cross method remains consistently the fastest for the considered range of tolerances.
For an even less smooth field with $\nu=1$, both the single-level and the multilevel QMC method become faster than the ALS-Cross method and the rate of ALS-Cross is closer to that of single-level QMC (Fig. \ref{fig:nu2}, left).
This indicates that it might be beneficial to use the ALS-Cross solver also in a multilevel framework.
A further investigation of this idea goes beyond the scope of this paper.
However, in terms of PDE solves (Fig.~\ref{fig:nu2}, right), ALS-Cross still hugely outperforms both the QMC approaches,
which may indicate potential even for rough coefficient fields for more complex 3D PDE problems.

\begin{figure}[t]
\centering
\caption{CPU time (left) and number of deterministic solves (right)
  vs.~$\log_{10}\widehat \eps_{\text{tot}}$ for the log-normal coefficient with $\nu=1$ (\#  solves in mlqmc computed according to Rem.~\ref{rem:norm_nsolve}).}
\label{fig:nu2}
\resizebox{\figscale\linewidth}{!}{
\begin{tikzpicture}
  \begin{axis}[%
  xmode=normal,
  ymode=log,
  ylabel=CPU time,
  legend style={at={(0.01,0.01)},anchor=south west},
  x filter/.code={\pgfmathparse{log10(\pgfmathresult*1)}\pgfmathresult},
  ]
   \addplot+[]  coordinates {
                              (6.8527e-03, 2.6134e+02)
                              (1.1868e-03, 1.3259e+04)
                              (2.5261e-04, 3.0602e+05)
                             };\addlegendentry{tt}
   \node[anchor=west] at (axis cs:-2.3,5e2) {\footnotesize\color{blue} 2.14};
   \addplot+[]  coordinates{
                              (3.5657e-03, 5.1410e+02)
                              (1.0235e-03, 9.3158e+03)
                              (2.5546e-04, 1.6698e+05)
                            };\addlegendentry{qmc}
   \node[anchor=east] at (axis cs:-2.55,8e2) {\footnotesize\color{red} 2.20};

   \pgfplotsset{cycle list shift=1}
   \addplot+[black]  coordinates {
                              (3.3828e-03, 1.0547e+03)
                              (8.0566e-04, 1.5361e+04)
                              (1.7126e-04, 1.7521e+05)
                             };\addlegendentry{mlqmc}
   \node[anchor=east] at (axis cs:-3.5,5e4) {\footnotesize\color{black} 1.71};
  \end{axis}
 \end{tikzpicture}
}
\resizebox{\figscale\linewidth}{!}{
 \begin{tikzpicture}
  \begin{axis}[%
  xmode=normal,
  ymode=log,
  ylabel=\# solves,
  x filter/.code={\pgfmathparse{log10(\pgfmathresult*1)}\pgfmathresult},
  ]
   \addplot+[]  coordinates{
                              (6.8527e-03, 2586)
                              (1.1868e-03, 5612)
                              (2.5261e-04, 11436)
                            };
   \addplot+[]  coordinates{
                              (3.5657e-03,  131072)
                              (1.0235e-03,  524288)
                              (2.5546e-04, 2097152)
                            };
   \pgfplotsset{cycle list shift=1}
   \addplot+[black] coordinates{
                              (3.3828e-03, 261888)
                              (8.0566e-04, 7.2059e+05)
                              (1.7126e-04, 1.2632e+06)
                                };
  \end{axis}
 \end{tikzpicture}
}
\end{figure}


\begin{figure}[t]
\centering
\caption{CPU time (left) and number of deterministic solves (right)
  vs.~$\log_{10}\widehat \eps_{\text{tot}}$ for log-normal coefficients with
  $\nu=3$, varying the correlation length $1/k_0$ (top) and the
  variance $\sigma^2$ (bottom). The number of solves in mlqmc is again
  computed according to Remark \ref{rem:norm_nsolve}.}
\label{fig:nu4_l_sigma}
\resizebox{\figscale\linewidth}{!}{
\begin{tikzpicture}
  \begin{axis}[%
  xmode=normal,
  ymode=log,
  ylabel=CPU time,
  xmin=-6.6,
  legend style={at={(0.01,0.99)},anchor=north west},
  x filter/.code={\pgfmathparse{log10(\pgfmathresult*1)}\pgfmathresult},
  ]
   \addplot+[black,solid,mark=*,mark options={black}] coordinates{ 
     (1.0278e-02,  4.9601e+01)
     (1.2505e-03,  5.7663e+01)
     (4.5851e-04,  3.3853e+02)
     (9.3987e-05,  1.7413e+03)
     (6.2240e-06,  8.9567e+03)
   }; \addlegendentry{$\substack{\mathtt{tt} \\ \ell=0.5}$}
   \addplot+[black,dashed,mark=*,mark options={black}] coordinates{  
     (6.1259e-03, 2.0987e+03)
     (1.3841e-03, 1.5750e+04)
     (2.2783e-04, 6.3816e+04)
     (1.5306e-04, 2.3088e+05)
   }; \addlegendentry{$\substack{\mathtt{mlqmc} \\ \ell=0.5}$}

   \addplot+[blue,solid,mark=square*,mark options={blue}] coordinates{  
     (7.1394e-03, 6.4855e+01)
     (2.8943e-03, 1.2618e+02)
     (4.5536e-04, 6.9657e+02)
     (1.6965e-04, 2.8585e+03)
   }; \addlegendentry{$\substack{\mathtt{tt} \\ \ell=0.3}$}
   \addplot+[blue,dashed,mark=square*,mark options={blue}] coordinates{ 
     (5.2568e-03, 1.0560e+03)
     (1.2528e-03, 1.4763e+04)
     (2.3902e-04, 7.0087e+04)
     (1.0568e-05, 2.1905e+05)
   }; \addlegendentry{$\substack{\mathtt{mlqmc} \\ \ell=0.3}$}

   \addplot+[orange,solid,mark=triangle*,mark options={orange}] coordinates{  
     (6.7005e-03, 5.2584e+02)
     (2.6785e-03, 1.4291e+03)
     (2.5120e-04, 4.8774e+03)
   }; \addlegendentry{$\substack{\mathtt{tt} \\ \ell=0.2}$}
   \addplot+[orange,dashed,mark=triangle*,mark options={orange}] coordinates{ 
     (5.5565e-03, 4.2082e+03)
     (1.1019e-03, 4.6329e+04)
     (4.5403e-05, 1.5520e+05)
   }; \addlegendentry{$\substack{\mathtt{mlqmc} \\ \ell=0.2}$}

   \addplot+[red,solid,mark=diamond*,mark options={red}] coordinates{  
     (9.0128e-03, 6.8319e+03)
     (1.2568e-03, 1.3061e+05)
   }; \addlegendentry{$\substack{\mathtt{tt} \\ \ell=0.1}$}
   \addplot+[red,dashed,mark=diamond*,mark options={red}] coordinates{ 
     (3.9637e-03, 4.2544e+03)
     (8.3636e-04, 3.4904e+04)
     (4.8917e-05, 2.8510e+05)
   }; \addlegendentry{$\substack{\mathtt{mlqmc} \\ \ell=0.1}$}
  \end{axis}
 \end{tikzpicture}
}
\resizebox{\figscale\linewidth}{!}{
 \begin{tikzpicture}
  \begin{axis}[%
  xmode=normal,
  ymode=log,
  ylabel=\# solves,
  xmin=-6.6,
  x filter/.code={\pgfmathparse{log10(\pgfmathresult*1)}\pgfmathresult},
  ]

   \addplot+[black,solid,mark=*,mark options={black}] coordinates{ 
     (1.0278e-02,   966         )
     (1.2505e-03,  1018         )
     (4.5851e-04,  1301         )
     (9.3987e-05,  1433         )
     (6.2240e-06,  1558         )
   };
   \addplot+[black,dashed,mark=*,mark options={black}] coordinates{  
     (6.1259e-03,    524224       )
     (1.3841e-03, 6.9755e+05      )
     (2.2783e-04, 5.2840e+05      )
     (1.5306e-04, 4.6702e+05      )
   };

   \addplot+[blue,solid,mark=square*,mark options={blue}] coordinates{  
     (7.1394e-03, 1650         )
     (2.8943e-03, 1656         )
     (4.5536e-04, 1832         )
     (1.6965e-04, 1878         )
   };
   \addplot+[blue,dashed,mark=square*,mark options={blue}] coordinates{ 
     (5.2568e-03,    262080       )
     (1.2528e-03, 6.6045e+05      )
     (2.3902e-04, 6.8337e+05      )
     (1.0568e-05, 4.6071e+05      )
   };

   \addplot+[orange,solid,mark=triangle*,mark options={orange}] coordinates{  
     (6.7005e-03, 3382         )
     (2.6785e-03, 3457         )
     (2.5120e-04, 3524         )
   };
   \addplot+[orange,dashed,mark=triangle*,mark options={orange}] coordinates{ 
     (5.5565e-03,   1048512       )
     (1.1019e-03, 2.0573e+06      )
     (4.5403e-05, 1.1521e+06      )
   };

   \addplot+[red,solid,mark=diamond*,mark options={red}] coordinates{  
     (9.0128e-03,  4624         )
     (1.2568e-03, 10270         )
   };
   \addplot+[red,dashed,mark=diamond*,mark options={red}] coordinates{ 
     (3.9637e-03,  1048512       )
     (8.3636e-04, 1.5465e+06     )
     (4.8917e-05, 2.8204e+06     )
   };

  \end{axis}
 \end{tikzpicture}
}\\
\resizebox{\figscale\linewidth}{!}{
\begin{tikzpicture}
  \begin{axis}[%
  xmode=normal,
  ymode=log,
  ylabel=CPU time,
  xmin=-6.2,ymin=3e0,
  legend style={at={(0.01,0.01)},anchor=south west},
  x filter/.code={\pgfmathparse{log10(\pgfmathresult*1)}\pgfmathresult},
  ]
   \addplot+[black,solid,mark=*,mark options={black}] coordinates{ 
             (4.2006e-03,   3.2653e+01)
             (7.0198e-04,   2.4518e+01)
             (4.2308e-04,   1.2175e+02)
             (3.9320e-05,   7.2275e+02)
             (1.6645e-05,   4.8143e+03)
   };\addlegendentry{$\substack{\mathtt{tt} \\ \sigma^2=0.3}$}
   \addplot+[black,dashed,mark=*,mark options={black}] coordinates{ 
           (7.7031e-04,   5.3567e+02)
           (2.8111e-04,   4.2130e+03)
           (5.2874e-05,   1.9929e+04)
           (1.3437e-05,   6.0118e+04)
           (5.0684e-06,   3.2593e+05)
   }; \addlegendentry{$\substack{\mathtt{mlqmc} \\ \sigma^2=0.3}$}

   \addplot+[blue,solid,mark=square*,mark options={blue}] coordinates{ 
               (6.1716e-03,   3.3268e+01)
               (1.0670e-03,   4.3222e+01)
               (4.9035e-04,   2.1017e+02)
               (5.5566e-05,   1.0962e+03)
               (6.1287e-06,   7.0860e+03)
   }; \addlegendentry{$\substack{\mathtt{tt} \\ \sigma^2=1}$}
   \addplot+[blue,dashed,mark=square*,mark options={blue}] coordinates{ 
           (1.4234e-03,   1.0526e+03)
           (3.5954e-04,   7.6645e+03)
           (8.6981e-05,   3.4941e+04)
           (2.1105e-05,   1.4681e+05)
   }; \addlegendentry{$\substack{\mathtt{mlqmc} \\ \sigma^2=1}$}

   \addplot+[red,solid,mark=triangle*,mark options={red}] coordinates{ 
              (5.2319e-03,  3.9501e+01)
              (1.7127e-03,  8.5206e+01)
              (4.2507e-04,  5.3724e+02)
              (1.3095e-04,  3.2496e+03)
   };\addlegendentry{$\substack{\mathtt{tt} \\ \sigma^2=3}$}
   \addplot+[red,dashed,mark=triangle*,mark options={red}] coordinates{ 
          (3.2174e-03,   2.1032e+03)
          (8.0115e-04,   1.3909e+04)
          (1.6433e-04,   5.9431e+04)
          (1.4872e-05,   2.0864e+05)
   }; \addlegendentry{$\substack{\mathtt{mlqmc} \\ \sigma^2=3}$}
  \end{axis}
 \end{tikzpicture}
}
\resizebox{\figscale\linewidth}{!}{
 \begin{tikzpicture}
  \begin{axis}[%
  xmode=normal,
  ymode=log,
  ylabel=\# solves,
  xmin=-6.2,
  x filter/.code={\pgfmathparse{log10(\pgfmathresult*1)}\pgfmathresult},
  ]
   \addplot+[black,solid,mark=*,mark options={black}] coordinates{ 
             (4.2006e-03,    337)
             (7.0198e-04,    475)
             (4.2308e-04,    647)
             (3.9320e-05,    900)
             (1.6645e-05,   1174)
   };
   \addplot+[black,dashed,mark=*,mark options={black}] coordinates{ 
           (7.7031e-04,      131008 )
           (2.8111e-04,   1.9066e+05)
           (5.2874e-05,   1.9660e+05)
           (1.3437e-05,   1.2850e+05)
           (5.0684e-06,   1.4516e+05)
   };

   \addplot+[blue,solid,mark=square*,mark options={blue}] coordinates{ 
               (6.1716e-03,    663)
               (1.0670e-03,    847)
               (4.9035e-04,   1047)
               (5.5566e-05,   1211)
               (6.1287e-06,   1414)
   };
   \addplot+[blue,dashed,mark=square*,mark options={blue}] coordinates{ 
           (1.4234e-03,     262080  )
           (3.5954e-04,   3.4725e+05)
           (8.6981e-05,   3.4278e+05)
           (2.1105e-05,   4.4552e+05)
   };

   \addplot+[red,solid,mark=triangle*,mark options={red}] coordinates{ 
              (5.2319e-03,  1269)
              (1.7127e-03,  1739)
              (4.2507e-04,  2152)
              (1.3095e-04,  2529)
   };
   \addplot+[red,dashed,mark=triangle*,mark options={red}] coordinates{ 
          (3.2174e-03,      524224 )
          (8.0115e-04,   6.2632e+05)
          (1.6433e-04,   7.9956e+05)
          (1.4872e-05,   6.1579e+05)
   };

  \end{axis}
 \end{tikzpicture}
}
\end{figure}
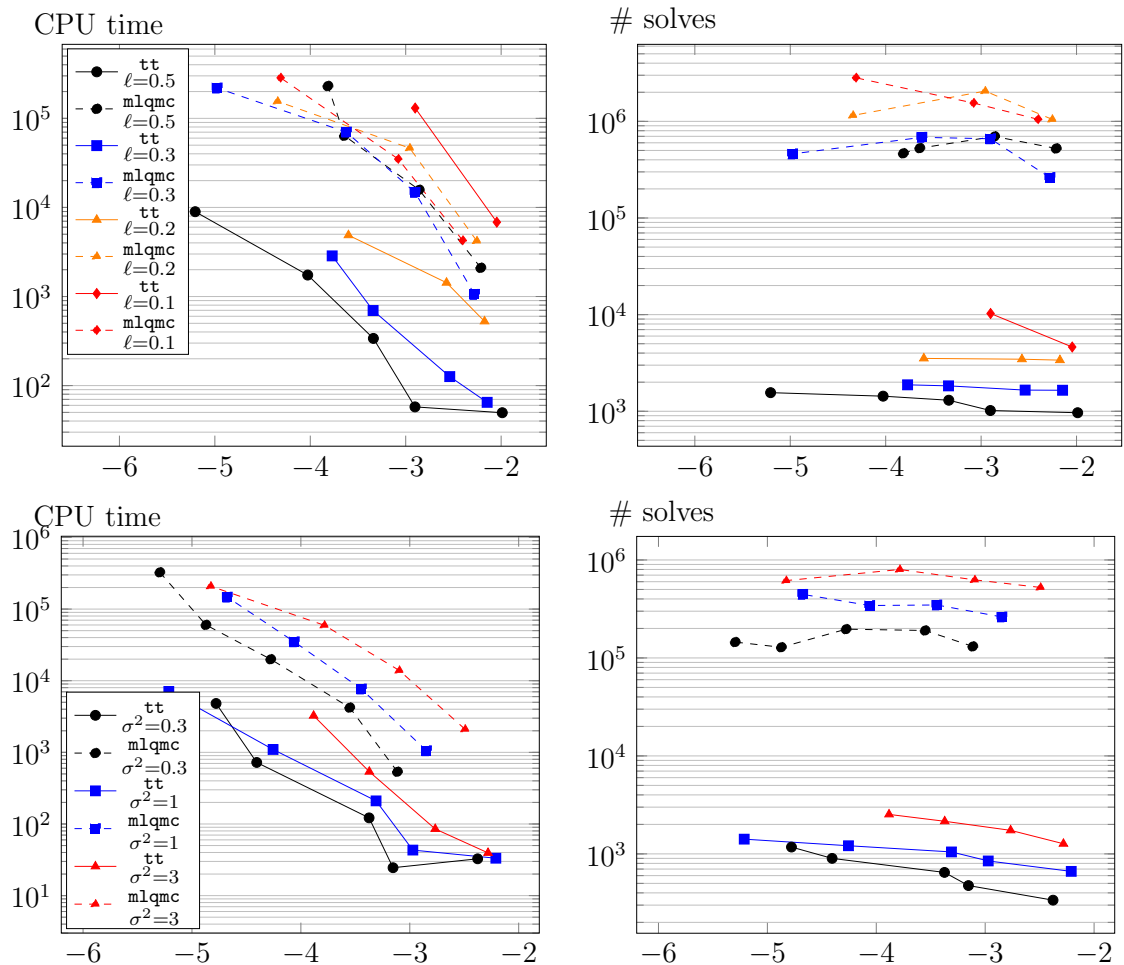
Finally, we investigate in Fig.~\ref{fig:nu4_l_sigma} the influence of
the parameters $k_0$ and $\sigma^2$ in our coefficient model \eqref{eq:kle_art} on the performance. These two parameters correspond to the
reciprocal of the correlation length and to the variance, respectively.
We consider only ALS-Cross and multilevel QMC here.
The results confirm what we expected: the TT ranks, and hence the complexity, are larger for a larger variance and for a smaller correlation length.
Consequently, the TT-based algorithm is more beneficial in the low variance, long correlation length regime.

\subsection{Estimation of probability densities}
\label{sec:num_PDF}

Finally, we compute probability density functions via the maximum
entropy method.
Let $\mathbf{Q}(\mathbf{y})$ be the QoI in \eqref{eq:QoI}, which we treat
as a random variable. Let us assume that $\mathbf{Q}$ has a probability density
function (PDF) $P(\mathbf{Q})>0$,
normalised such that the zeroth moment $\int_{\mathbb{R}} P(\mathbf{Q}) d\mathbf{Q} = 1$.
The $p$-th moment of $\mathbf{Q}$, for $p \ge 1$,
is then given by $\int_{\mathbb{R}} \mathbf{Q}^p P(\mathbf{Q}) d\mathbf{Q}$.
On the other hand, we have just described how to compute
approximations $Q_p$ for these moments from the parametric solution
\eqref{eq:QoI}.

To approximate the entire PDF of $\mathbf{Q}$ we can use the
maximum entropy method described in \cite{Kaverhad-maxentr-1986}.
It maximizes the Shannon
entropy $\mathbb{E}(-\ln(P(\mathbf{Q}))) = - \int_{\mathbb{R}}
P(\mathbf{Q}) \ln P(\mathbf{Q}) d\mathbf{Q}$, subject to matching the first
$S+1$ moments of $\mathbf{Q}$ (including the zeroth moment) to the
computed moments. By enforcing these constraints via Lagrange multipliers
$\lambda=(\lambda_0,\ldots,\lambda_S)$, the first-order criticality
conditions lead to an approximate PDF of the form
$P_S(\mathbf{Q}) = \exp\left(\sum_{p=0}^{S} \mathbf{Q}^p
  \lambda_p\right)$ where the Lagrange multipliers $\lambda$ satisfy
the $(S+1)$-dimensional nonlinear system
$$
\int_{\mathbb{R}} \mathbf{Q}^p P_S(\mathbf{Q}) d\mathbf{Q} = \left\{\begin{array}{ll}1, & p=0, \\ Q_p, & p=1,\ldots,S, \end{array}\right.
$$
which can be solved via an efficient Newton algorithm
\cite{Kaverhad-maxentr-1986}. So the main computational cost is in
computing the moments $Q_p$, as described above.

\begin{figure}[t]
\centering
\caption{PDFs and errors for different choices of $S$, the number of
  moments, and for log-normal
  (top) and log-uniform (bottom) coefficients with $\nu=3$ (left) and
  $\nu=1.5$ (right).}  
\label{fig:pdf}
\resizebox{\figscale\linewidth}{!}{
 \begin{tikzpicture}
  \begin{axis}[%
  ymode=normal,
  xmode=normal,
  xmin=0,xmax=0.4,
  ymin=0,ymax=13,
  ylabel={$P_S$},
  ytick pos=left,
  legend style={at={(0.01,0.2)},anchor=south west},
  ]
   \addplot+[no marks] table[header=true, x=t, y=f2]{pdf_nu4.dat}; \addlegendentry{\footnotesize $S=2$};
   \addplot+[no marks] table[header=true, x=t, y=f4]{pdf_nu4.dat}; \addlegendentry{\footnotesize $S=4$};
   \addplot+[no marks] table[header=true, x=t, y=f6]{pdf_nu4.dat}; \addlegendentry{\footnotesize $S=6$};
   \addplot+[no marks] table[header=true, x=t, y=f8]{pdf_nu4.dat}; \addlegendentry{\footnotesize $S=8$};
  \end{axis}
  \begin{axis}[%
   xmode=normal,
   ymode=log,
   ymin=8e-4,ymax=1.1e-1,
   axis y line*=right, y label style={at={(1.3,1.0)},anchor=south east,color=cyan!50!black},
   axis x line=none,
   ylabel={$\|P_S-P_8\|_2/\|P_8\|_2$},
   yminorgrids=false,
   nodes near coords, nodes near coords align={horizontal},point meta=explicit symbolic,
   yticklabel style={color=cyan!50!black},
   ]
   \addplot+[cyan!50!black,mark options={cyan!50!black}] coordinates{(2,3.7202e-02)[2] (3,2.2706e-02)[3] (4,4.3177e-03)[4] (5,1.8954e-03)[5] (6,1.8673e-03)[6] (7,1.6429e-03)[7]};
  \end{axis}
 \end{tikzpicture}
}
\resizebox{\figscale\linewidth}{!}{
 \begin{tikzpicture}
  \begin{axis}[%
  ymode=normal,
  xmode=normal,
  xmin=0,xmax=0.4,
  ymin=0,ymax=13,
  ylabel={$P_S$},
  ytick pos=left,
  legend style={at={(0.01,0.3)},anchor=south west},
  ]
   \pgfplotsset{cycle list shift=1};
   \addplot+[no marks] table[header=true, x=t, y=f4]{pdf_nu2.5.dat}; \addlegendentry{\footnotesize $S=4$};
   \addplot+[no marks] table[header=true, x=t, y=f6]{pdf_nu2.5.dat}; \addlegendentry{\footnotesize $S=6$};
   \addplot+[no marks] table[header=true, x=t, y=f8]{pdf_nu2.5.dat}; \addlegendentry{\footnotesize $S=8$};
   \addplot+[no marks,green!60!black] table[header=true, x=t, y=f10]{pdf_nu2.5.dat}; \addlegendentry{\footnotesize $S=10$};
  \end{axis}
  \begin{axis}[%
   xmode=normal,
   ymode=log,
   ymin=8e-4,ymax=1.1e-1,
   axis y line*=right, y label style={at={(1.3,1.0)},anchor=south east,color=cyan!50!black},
   axis x line=none,
   ylabel={$\|P_S-P_{10}\|_2/\|P_{10}\|_2$},
   yminorgrids=false,
   nodes near coords, nodes near coords align={horizontal},point meta=explicit symbolic,
   yticklabel style={color=cyan!50!black},
   ]
   \addplot+[cyan!50!black,mark options={cyan!50!black}] coordinates{(2,6.6721e-02)[2] (3,4.6422e-02)[3] (4,1.2227e-02)[4] (5,6.9211e-03)[5] (6,3.5513e-03)[6] (7,2.9258e-03)[7] (8,2.4723e-03)[8] (9,1.2829e-03)[9]};
  \end{axis}
 \end{tikzpicture}
}\\
\resizebox{\figscale\linewidth}{!}{
 \begin{tikzpicture}
  \begin{axis}[%
  ymode=normal,
  xmode=normal,
  xmin=0,xmax=0.4,
  ymin=0,ymax=11,
  ylabel={$P_S$},
  ytick pos=left,
  legend style={at={(0.99,0.99)},anchor=north east},
  ]
   \addplot+[no marks] table[header=true, x=t, y=f5]{pdf_nu4_u.dat}; \addlegendentry{\footnotesize $S=5$};
   \addplot+[no marks] table[header=true, x=t, y=f7]{pdf_nu4_u.dat}; \addlegendentry{\footnotesize $S=7$};
   \addplot+[no marks,green!50!black] table[header=true, x=t, y=f9]{pdf_nu4_u.dat}; \addlegendentry{\footnotesize $S=9$};
   \addplot+[no marks] table[header=true, x=t, y=f10]{pdf_nu4_u.dat}; \addlegendentry{\footnotesize $S=10$};
  \end{axis}
  \begin{axis}[%
   xmode=normal,
   ymode=log,
   ymin=3e-3,ymax=0.8,
   axis y line*=right, y label style={at={(1.3,1.0)},anchor=south east,color=cyan!50!black},
   axis x line=none,
   ylabel={$\|P_S-P_{10}\|_2/\|P_{10}\|_2$},
   yminorgrids=false,
   nodes near coords, nodes near coords align={horizontal},point meta=explicit symbolic,
   yticklabel style={color=cyan!50!black},
   ]
   \addplot+[cyan!50!black,mark options={cyan!50!black}] coordinates{(5, 1.5211e-01)[5] (6, 5.6133e-02
)[6] (7, 5.3032e-02)[7] (8, 1.7707e-02)[8] (9, 6.4909e-03)[9]};
  \end{axis}
 \end{tikzpicture}
}
\resizebox{\figscale\linewidth}{!}{
 \begin{tikzpicture}
  \begin{axis}[%
  ymode=normal,
  xmode=normal,
  xmin=0,xmax=0.4,
  ymin=0,ymax=11,
  ylabel={$P_S$},
  ytick pos=left,
  legend style={at={(0.99,0.99)},anchor=north east},
  ]
   \addplot+[no marks] table[header=true, x=t, y=f5]{pdf_nu2.5_u.dat}; \addlegendentry{\footnotesize $S=5$};
   \addplot+[no marks] table[header=true, x=t, y=f7]{pdf_nu2.5_u.dat}; \addlegendentry{\footnotesize $S=7$};
   \addplot+[no marks,green!50!black] table[header=true, x=t, y=f9]{pdf_nu2.5_u.dat}; \addlegendentry{\footnotesize $S=9$};
   \addplot+[no marks] table[header=true, x=t, y=f10]{pdf_nu2.5_u.dat}; \addlegendentry{\footnotesize $S=10$};
  \end{axis}
  \begin{axis}[%
   xmode=normal,
   ymode=log,
   ymin=3e-3,ymax=0.8,
   axis y line*=right, y label style={at={(1.3,1.0)},anchor=south east,color=cyan!50!black},
   axis x line=none,
   ylabel={$\|P_S-P_{10}\|_2/\|P_{10}\|_2$},
   yminorgrids=false,
   nodes near coords, nodes near coords align={horizontal},point meta=explicit symbolic,
   yticklabel style={color=cyan!50!black},
   ]
   \addplot+[cyan!50!black,mark options={cyan!50!black}] coordinates{(5, 9.2573e-02)[5] (6, 3.7990e-02
)[6] (7, 1.5096e-02)[7] (8, 5.4980e-03)[8] (9, 4.0285e-03)[9]};
  \end{axis}
 \end{tikzpicture}
}
\end{figure}
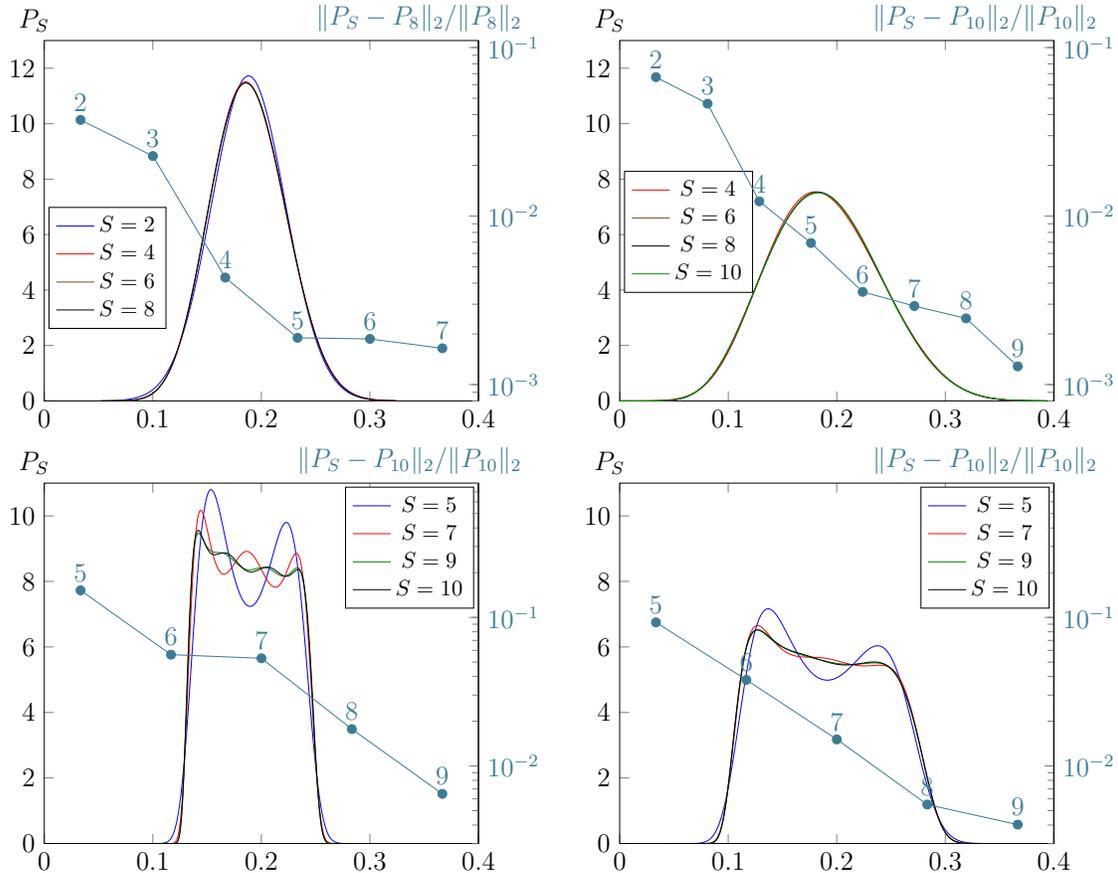

In Fig.~\ref{fig:pdf} we investigate two values of the smoothness parameter $\nu$ for the log-normal and log-uniform cases.
The left axes in Fig.~\ref{fig:pdf} show the PDFs for different numbers of moments.
The right axes show the relative error in the PDF approximation $P_S$
with $S$ moments. As the exact solution in those error estimates, we
use a computation with a larger number of moments ($8$ for the top
left figure, $10$ for the other plots).
We see that up to $10$ moments might be necessary for approximating
the PDF to achieve a relative error below $10^{-2}$.
For a smaller relative error, more moments computed with a smaller tolerance
$\delta$ would be necessary. This is the rationale behind
computing the full solution in the TT format, particularly for more
complex PDEs, since the number of deterministic problems the TT method needs to solve is several orders of magnitude lower than in the other methods, for all considered parameters.

In most of our numerical experiments, the TT method is faster than multilevel QMC for lower accuracies,
since the TT approach needs much fewer deterministic solves.
On the other hand, in these experiments the MLQMC cost is dominated by the
work on the coarsest level,
while the finer level corrections require less than a thousand samples.
Thus, the multilevel correction idea could provide a more efficient way to
enhance the accuracy of the TT solution, especially for rougher
coefficients (i.e. for smaller values of $\nu$), similar to
\cite{nobile-mlmcwithcv-2015}.
More importantly for applications, the TT solution of the forward
problem can also be used as an efficient and accurate surrogate model
for sampling-based Bayesian approaches to high-dimensional inverse
problems, reducing the number of expensive PDE solves. Both these
extensions will be the focus of our future research.

\bibliographystyle{siam}
\bibliography{stochastic,tensor,our,iter,dmrg,misc}

\end{document}